\documentclass[twocolumn]{autart} 
\usepackage{graphicx,float,wrapfig}
\usepackage{epstopdf}
\usepackage{epsfig}
\usepackage{amssymb,amsmath,amsfonts,layout}
\usepackage{multirow}
\usepackage{subfigure}
\usepackage{url}
\usepackage{enumerate}
\usepackage{cite}

\begin{document}

\begin{frontmatter}

\title{Distributed Control of Inverter-Based Lossy Microgrids for Power Sharing and Frequency Regulation Under Voltage Constraints} 


\author[USA]{Chin-Yao Chang}\ead{chang.981@osu.edu} and    
\author[USA]{Wei Zhang}\ead{zhang@ece.osu.edu}   

\address[USA]{Department of Electrical and Computer Engineering, Ohio State University, Columbus, OH 43210, USA}

\begin{keyword}                           
Microgrid Control, Droop Control, Frequency Synchronization, Power Sharing, Voltage Regulation          
\end{keyword}     


\begin{abstract}
This paper presents a new distributed control framework to coordinate inverter-interfaced distributed energy resources (DERs) in island microgrids. We show that under bounded load uncertainties, the proposed control method can steer the microgrid to a desired steady state with synchronized inverter frequency across the network and proportional sharing of both active and reactive powers among the inverters. We also show that such convergence can be achieved while respecting constraints on voltage magnitude and branch angle differences. The controller is robust under various contingency scenarios, including loss of communication links and failures of DERs. The proposed controller is applicable to lossy mesh microgrids with heterogeneous R/X distribution lines and reasonable parameter variations. Simulations based on various microgrid operation scenarios are also provided to show the effectiveness of the proposed control method.
\end{abstract}

\end{frontmatter}

\section{Introduction}
Microgrids are low voltage power networks comprised of distributed generations (DGs), energy storages systems (ESSs), and loads that can operate in either grid-connected or island mode. Distributed generation contributes on-site and clean energy, which is expected to make power networks more robust, efficient and environmentally friendly~\cite{pepermans2005distributed,ackermann2001distributed}. Energy storage systems are considered as an important resource to benefit the power networks by smoothing real time imbalance between generation and demand~\cite{yang2011electrochemical}. Some storage devices such as freewheel and battery packs can be integrated with intermittent DGs to regulate the power injection to a power network~\cite{prasad2006optimization,vazquez2010energy}. Demand side appliances such as plug-in hybrid electric vehicles (PHEV) and thermostatically controlled loads (TCLs) can also be viewed as energy storage resources. Those ``storage'' appliances can be coordinated to provide ancillary services to the main grid~\cite{zhao2013optimal, zhao2014design, zhang2013aggregated}. The proximity of DGs and ESSs to loads in a microgrid allows for a transition to the island mode during faults on the main grid. Such a transition may also be triggered by efficiency or reliability incentives, (see \cite{nourai2010batteries, piagi2006autonomous}).

Distributed energy resources (DERs) such as DGs and ESSs connect to the microgrid through DC/AC or AC/AC inverters. During the island mode, the inverters are typically operated as voltage source inverters (VSIs). These VSIs need to be controlled cooperatively to achieve desired performance and reliability properties. In AC networks, voltage magnitude and angle difference between connected buses should be regulated in some bounded ranges for system security and stability. Frequency synchronization to a nominal value is also crucial for grid connection and stability purposes. Besides frequency and voltage regulation, sharing of active and reactive power is also considered as important control objectives in microgrids~\cite{mohamed2008adaptive, schiffer14_dvc}. They require that the power injection into the microgrid from DERs is proportional to the nominal value defined by economics or other incentives, while satisfying load demands~\cite{schiffera2013conditions}. Power sharing enables effective utilization of limited generation resources and prevents overloading~\cite{katiraei2008microgrids}.

To achieve the aforementioned objectives, a microgrid is typically controlled using a hierarchal structure including primary, secondary, and tertiary controls \cite{guerrero2011hierarchical, Bidram2012, guerrero2013advanced, dorfler2014breaking}, which is similar to the one used in the traditional power systems. The primary droop control of a microgrid maintains the voltage and frequency stability while balancing the generation and load with proper power sharing. The secondary controller compensates the voltage and frequency deviations from their reference values. The tertiary control establishes the optimal power sharing between inverters in both islanding and grid-connected modes. 

The primary droop is generally a decentralized controller that adjusts the voltage frequency and magnitude of each inverter in response to active and reactive power deviations from their nominal values. Various droop methods are proposed to achieve proportional active and reactive load power sharing \cite{mohamed2008adaptive, yao2011design, tabatabaee2011investigation, iyer2010generalized, vasquez2013modeling, ahn2010power, lee2013new}. However, this is often achieved at the cost of sacrificing other control objectives such as voltage and frequency regulation. The secondary control utilizes either centralized or decentralized communication infrastructures to restore frequency and voltage deviation induced by the primary droop. Most of the existing secondary control methods require centralized communications \cite{li2009accurate, he2012enhanced, mehrizi2010potential}. On the other hand, decentralized secondary control has recently been proposed to avoid single point of failure \cite{shafiee2014distributed}. The combined operations of the primary and secondary control require separation of time scale, resulting in slow dynamics that cannot effectively handle fast-varying loads \cite{simpson2013synchronization}. In addition, the secondary control may destroy the proportional power sharing established in the primary control layer \cite{simpson2013stability}. One possible solution is to adopt distributed or decentralized control structure for primary and secondary control layers to improve performance and support plug-and-play operation of the microgrid~\cite{dorfler2014breaking}.

Many existing primary and secondary control methods rely on small signal linearization for stability analysis, which is vulnerable to parameter variations and change of operating points. Only several recent works \cite{simpson2013stability, schiffera2013conditions, ainsworth2013structure} have rigorously analyzed the stability of microgrid with droop-controlled inverters. In particular, \cite{simpson2013stability} derives a necessary and sufficient condition for the stability under primary droop control. The authors have also proposed a distributed averaging controller to fix the time scale separation issue between the primary and secondary control layers. In \cite{schiffera2013conditions} and \cite{ainsworth2013structure}, stability conditions of lossless mesh microgrids have been provided. Despite their advantages, these nonlinear methods still suffer from several common limitations. First, all the nonlinear analyses mentioned above only focus on lossless microgrids with purely inductive distribution lines. The results may not be applicable for microgrids with heterogeneous and mixed R/X ratio lines, which is common in low voltage microgrids~\cite{li2009accurate}. Secondly, since only frequency droop is carefully analyzed, reactive power sharing is often not guaranteed. 

To address the aforementioned limitations of the existing works, we propose a distributed control framework to coordinate VSIs in an island AC microgrid. The proposed control adjusts each inverter frequency and voltage magnitude based on the active/reactive power measurements of its neighbors. We first show that the particular control structure ensures that any equilibrium of the closed-loop system results in the desired power sharing and frequency synchronization. Secondly, conditions for power sharing and synchronized frequency respecting voltage constraints are provided. The proposed controller can be applied to both radial and mesh microgrids with mixed R/X ratios. Furthermore, the proposed controller requires no separation of time scale and can tolerate reasonable parameter variations. To the authors' knowledge, most existing control framework cannot achieves active/reactive power sharing while respecting voltage and frequency regulation for a mesh micogrid with mixed R/X ratio lines.

To demonstrate the robustness of the proposed distributed controller, we also study the control performance under partial communication failures and the plug-and-play operations. We will show that as long as the communication network remains connected, all the desired properties including power sharing and frequency and voltage regulation still hold in these contingency scenarios. This effectively demonstrates the robustness of the proposed distributed controller. It is worth to mention that the proposed framework may require faster communications among the VSIs than the traditional secondary control. However, such communication requirement is reasonable for most microgrid control systems \cite{xin2011self, gungor2011smart, laaksonen2010protection}. 

The rest of this paper is organized as follows. Section~\ref{sec: ProbSet} formulates the microgrid control problem. Sufficient conditions for the solvability of the proportional power sharing problem respecting voltage constraints are also provided. The proposed distributed control framework is developed in Section \ref{sec:ContrDesign}. Robustness of the distributed controller under loss of  communication links or failures of DERs is studied in Section \ref{sec:FlexOperMicro}. In Section \ref{sec:SimResults}, we validate the proposed controller through simulations under various microgrid operating scenarios, including abrupt changes of loads and loss of one VSI. Some concluding remarks are given in Section \ref{sec: conclusion}.

\textbf{Notation}
Define $\mathbb{R}_+$ and $\mathbb{R}_-$ as positive and negative real numbers, respectively. Denote $[n]:= \{1,2,...,n\}$. Given a set $\mathcal{V}$, let $|\mathcal{V}|$ and $2^{\mathcal{V}}$ be its cardinality and power set, respectively. Denote the diagonal matrix of a vector $x$ as $diag(x)$. For a set of vectors $x_i$, $i \in \mathcal{I}$, let $\{x_i, i \in \mathcal{I} \}$ be the augmented vector of $x_i$ collecting all $i \in \mathcal{I}$. Given a polyhedron $\mathcal{B} \in \mathbb{R}^n$, let $v(\mathcal{B})$ be the vertex set of $\mathcal{B}$. For a closed set $F \subseteq \mathbb{R}^n$, int($F$) and $\partial F$ are the interior and the boundary of $F$. The distance between a point $f \in \mathbb{R}^n$ and the set $F$ is denoted as $d(f,F):= \text{inf}\{||f - \bar{f}||_2\; | \bar{f} \in F \}$. Define ${\bf{1}}_n \in \mathbb{R}^n $ and ${\bf{0}}_n \in \mathbb{R}^n$ as the vectors with all the elements being ones and zeros, respectively. For a symmetric matrix $A$, let $\lambda(A)$ and $\underline{\lambda}(A)$ be the spectrum and minimal eigenvalue of $A$, respectively. Denote $A \otimes B$ as the tensor product between matrices $A$ and $B$. Let null($A$) be the null space a matrix $A$. 

\section{Problem Formulation}
\label{sec: ProbSet}
\begin{figure}
\begin{center}
\centerline{\scalebox{0.7}{\includegraphics{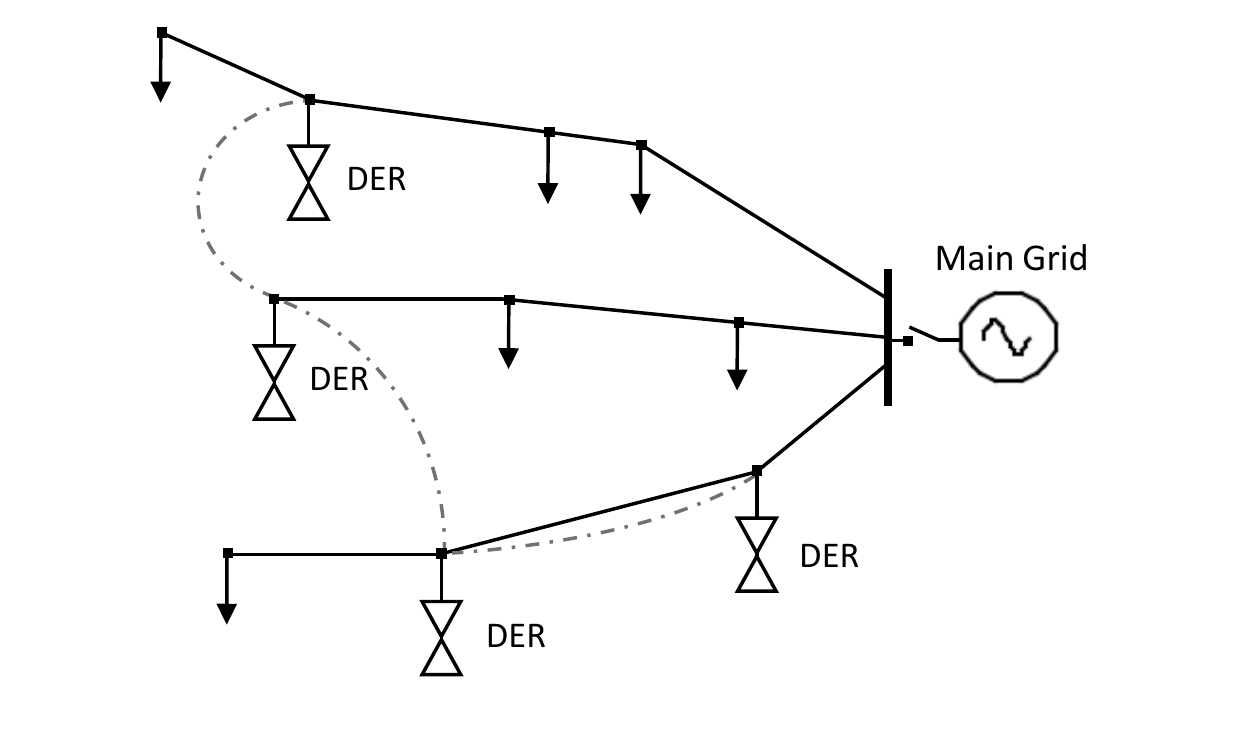}}}
\caption{A connected microgrid network. The dash lines represent the communication links and the solid lines represent the distribution lines connecting VSIs, loads, and the main grid.}
\label{fig:MicorgridTopology}
\end{center}
\end{figure}

In this paper, we consider a connected island microgrid network as shown in Fig. \ref{fig:MicorgridTopology}. An island microgrid is represented by a connected and undirected graph $\mathcal{G} = (\mathcal{V},\mathcal{E})$, where $\mathcal{V}$ is the set of buses (nodes) and $\mathcal{E} \subseteq \mathcal{V}\times \mathcal{V}$ is the set of distribution lines (edges) connecting the buses. The set of buses is partitioned into two parts, inverter buses $\mathcal{V}_I$ and load buses $\mathcal{V}_L$. Let $n_I = |\mathcal{V}_I|$, $n_L = |\mathcal{V}_L|$ and $n = |\mathcal{V}|$. The magnitude and phase angle of the bus voltage are denoted as $E_i$ and $\theta_i$, respectively. Let $x_i \triangleq [\theta_i, E_i]^T$ be the state vector at bus $i$, and let $x_I \triangleq \{x_i, i \in \mathcal{V}_I \}$ and $x_L = \{x_i, i \in \mathcal{V}_L \}$ be the inverter bus state vector and load bus state vector, respectively. The overall system state vector is denoted by $x = [x_I^T, x_L^T]^T$ and will be referred to as the system {\em voltage profile}. 
For each bus $i\in \mathcal{V}$, let $P_i$ and $Q_i$ be the active and reactive power injections at bus $i$. Given the admittance matrix $\textbf{Y}~\in~\mathbb{C}^{n \times n}$ of the microgrid, the active and reactive power injections are related to the voltage profile $x$ by the power flow equations \cite{glover2011power}
\begin{align} 
\label{eq:PFEqua}
\left\{
\begin{array}{ll}
P_{i}(x) &= E_i \sum_{j \in \mathcal{V}} Y_{ij}E_j \text{cos}(\theta_i - \theta_j - \phi_{ij})  \\ 
Q_{i}(x) &= E_i \sum_{j \in \mathcal{V}} Y_{ij}E_j \text{sin}(\theta_i - \theta_j - \phi_{ij}), 
\end{array}
\right.
\end{align}
where $Y_{ij} = ||\textbf{Y}_{ij}||_2$ and $\phi_{ij} = \angle \textbf{Y}_{ij}$ are the magnitude and the phase angle of the admittance matrix element $\textbf{Y}_{ij}$. 

We distinguish the voltage at inverter and load buses in our formulation due to their different characteristics. For inverter buses, there are standard methods to control the voltage magnitude and frequency (\cite{kim2010feedback}, \cite{yang2011robust}). Typically, these methods can track a given inverter voltage reference almost instantaneously. Therefore, an inverter is often modeled as a controlled voltage source behind a reactance~\cite{zhong2012control}. We also adopt such a model and consider $x_I$ can be fully controlled. In contrast to $x_I$, voltage at load buses $x_L$ is uncertain. The voltage dynamics of the load buses are assumed to satisfy the following condition
\begin{assum}
\label{assum:Bddxldot}
$|| \dot{x}_L ||_2 \leq \kappa || \dot{x}_I ||_2$,
for some $\kappa > 0$ 
\end{assum}
where $\kappa >0$ is a constant determined by load properties and microgrid topology. More explanation about this assumption is provided in Appendix~\ref{eq:Bddx_L}. Throughout this work, we will focus on constant power or constant impedance loads so that Assumption~\ref{assum:Bddxldot} holds. 
\begin{rem}
\normalfont{Under multi-agent or centralized control framework, load voltage is often assumed to be measurable and known during controller design~\cite{colson2011algorithms, corsi2008real}. In this paper, we consider a more general scenario, where load voltage is viewed as an unknown variables with only Assumption~\ref{assum:Bddxldot} being involved in the controller design. }
\end{rem}

Given nominal active and reactive power injections $P_i^*$ and $Q_i^*$, $i \in \mathcal{V}_I$, it is desired that the power injection of the inverters share uncertain loads proportionally to their nominal value:
\begin{defn}{\textnormal{\textbf{(Proportional Power Sharing)}}}
The active and reactive power are proportionally shared among the buses $j \in \mathcal{V}_I$ if 
\begin{align}
\label{eq:AcPSharingPre}
\frac{{P}_{j}(x)}{P^*_j} = \frac{{P}_{k}(x)}{P^*_k}, \;\;
\frac{{Q}_{j}(x)}{Q_{j}^*} = \frac{{Q}_{k}(x)}{Q_{k}^*}, \; j,k \in \mathcal{V}_I.
\end{align}
\end{defn}
The power sharing condition (\ref{eq:AcPSharingPre}) imposes a constraint on the system voltage profile. We define this constraint set as
\begin{align}  
\label{eq:AcPSharing}
x \in \mathcal{X}_S:= \{x | \text{Eq. (\ref{eq:AcPSharingPre}) holds} \}.
\end{align}
In addition to condition (\ref{eq:AcPSharing}), the control and operation of a microgrid has to respect its branch angle difference limits and voltage magnitude constraint. The branch angle difference between all connected buses is typically required to be bounded by a given constant $\gamma \in [0, \frac{\pi}{2})$. The upper bound $\gamma$ is derived based on the maximum current allowable on each distribution line $I_{jk}^{\max}$ (see~\cite{wu1982steady} and \cite{zimmerman2011matpower}). In addition to branch angle difference, the voltage magnitude also needs to stay inside some secure operation range~\cite{glover2011power}. Denote $[\underline{E}_i,\bar{E}_i]$ as a given desired range of the voltage magnitude of bus $i$. Both branch angle difference and voltage magnitude requirements impose a constraint on the system voltage profile defined below:
\begin{align*}
& \mathcal{X}_{\Theta}= \{x | |\theta_i - \theta_j| \leq \gamma, \; \forall \{i,j\} \in \mathcal{E} \}, \\
& \mathcal{X}_E= \{x | \underline{E}_i \leq E_i \leq \bar{E}_i \; \forall{i \in \mathcal{V}} \},
\end{align*}
These two constraints will be referred to as the security constraint of the microgrid
\begin{defn}{\textnormal{\textbf{(Security Constraints)}}}
We say that a microgrid satisfies the security constraints if
\begin{align}
\label{eq:SecConst}
& x \in \mathcal{X}_c:=  \mathcal{X}_{E} \cap \mathcal{X}_{\Theta}.
\end{align}
\end{defn}
In addition to proportional power sharing and voltage regulation, another important microgrid control objective is known as frequency regulation. Frequency regulation is defined as synchronization without deviations from the nominal value, specifically, $\dot{\theta}_i = \omega_0, \; \forall i \in \mathcal{V}$, where $\omega_0$ is a predefined nominal frequency of the microgrid. Note that the power flow equations (\ref{eq:PFEqua}) and the security constraints are invariant with respect to rigid rotation of ${\theta_i}$ of all buses. We can select a reference frame rotating at angular frequency $\omega_0$ while preserving all properties in Eq. (\ref{eq:PFEqua})-(\ref{eq:SecConst}). With the rotating reference, the frequency regulation condition is reduced to
\begin{align}
\label{eq:FReg}
\dot{\theta}_i = 0, \; \forall i \in \mathcal{V}.
\end{align}

As discussed in~\cite{guerrero2013advanced}, the controller design to meet the requirements specified in Eq. (\ref{eq:AcPSharing})-(\ref{eq:FReg}) inevitably requires communication networks. In this paper, we employ a distributed communication structure similar to~\cite{schiffer14_dvc} and \cite{simpson2013synchronization}, where each inverter can communicate with its neighboring inverters to share its local measurements as shown in Fig. \ref{fig:MicorgridTopology}. Let $\mathcal{G}_c = (\mathcal{V}_I, \mathcal{E}_c)$ be a connected simple graph of the communication network, where each edge $\{i,j\} \in \mathcal{E}_c$ represents an available communication link between buses $i$ and $j$. Let $\mathcal{N}_i := \{j| \{i,j\} \in \mathcal{E}_c \}\cup\{i\}$ be the set of neighbors of bus $i$, (including bus $i$ itself). An inverter $i$ has access to the measurements at every inverter bus $j \in \mathcal{N}_i$, including $P_j$ and $Q_j$. 

Since each inverter is modeled as a VSI, the microgrid-level coordination control for each inverter $i$ reduces to the determination of appropriate voltage frequency and magnitude setpoints. The actual frequency $\dot\theta_i$ and magnitude $E_i$ can track these setpoints almost instantaneously. The challenge here lies in that the load is uncertain and different load conditions require different voltage profile $x_I$ in order to satisfy constraints (\ref{eq:AcPSharing})-(\ref{eq:FReg}). Define $S_i = [{P}_i/P_{i}^*, {Q}_i/Q_{i}^*]^T$, $S_{\mathcal{N}_i}=\{ S_j, j \in \mathcal{N}_i \}$ and $E_{\mathcal{N}_i}=\{ E_j, j \in \mathcal{N}_i \}$. Our goal is thus to design a controller for each inverter $i$ that can automatically find the desired voltage vector $x_i$ based on local information $(S_{\mathcal{N}_i},E_{\mathcal{N}_i})$. Towards this end, we propose to dynamically update $x_i$ as follows
\begin{align}
\label{eq:PreController}
\left\{
\begin{array}{ll}
& \dot{x}_i = \mu_i(S_{\mathcal{N}_i}(x),E_{\mathcal{N}_i}), \; \forall i \in \mathcal{V}_I \\ 
& \text{subj. to } x(t) \in \mathcal{X}_c, \; \forall t\geq 0,  
\end{array}
\right.
\end{align}
where $\mu_i$ is the control law of inverter $i$ to be designed. Note that the above control structure corresponds to directly assigning frequency $\dot\theta_i$ based on local information, while dynamically updating voltage magnitude $E_i$ through simple integrator dynamics. Such structure is commonly used in the literature of microgrid control, (see \cite{schiffer14_dvc} and \cite{simpson2013voltage}). The constraint $x(t) \in \mathcal{X}_c$ is imposed to ensure that the security constraints are always satisfied. 

Define $\mathcal{X}_e:= \mathcal{X}_c \cap \mathcal{X}_S$, then conditions (\ref{eq:AcPSharing}) and (\ref{eq:SecConst}) hold when $x \in \mathcal{X}_e$. Under Assumption~\ref{assum:Bddxldot}, condition (\ref{eq:FReg}) holds when $\dot{x}_i = 0$ for all $i \in\mathcal{V}_I$. Our goal becomes to designing $\mu_i$ such that $\mathcal{X}_e$ forms an equilibrium set of system (\ref{eq:PreController}). In addition, we also want to achieve an exponential convergence to $\mathcal{X}_e$ for some initial $x(0)$. If such a controller is found, it can steer the microgrid to the desired steady state where conditions (\ref{eq:AcPSharing})-(\ref{eq:FReg}) hold. In the rest of this paper, we will first develop the control law $\mu_i$ such that $\mathcal{X}_e$ forms the equilibrium set of the system (\ref{eq:PreController}), and then derive conditions to ensure the exponential convergence of $\mathcal{X}_e$.

\section{A Distributed Microgrid Control Framework}
\label{sec:ContrDesign}

In this section, we propose a distributed control framework to coordinate the inverters in an island AC microgrid to accomplish the control objectives (\ref{eq:AcPSharing})-(\ref{eq:FReg}). We first provide sufficient conditions to ensure $\mathcal{X}_e\neq \emptyset$. A control design framework is then developed.

\subsection{Existence of Solutions}
A minimum requirement for the controller design is the existence of a voltage profile $x$ satisfying all the constraints, i.e., $\mathcal{X}_e\neq\emptyset$. Existing methods in the literature often directly assume this condition holds~(\cite{tabatabaee2011investigation}, \cite{vasquez2013modeling}). Here, we provide a brief discussion and a set of sufficient conditions to guarantee the non-emptiness of $\mathcal{X}_e$. The existence of the voltage profile $x$ satisfying conditions (\ref{eq:AcPSharing}) and (\ref{eq:SecConst}) involves solving nonlinear algebraic power flow equations (\ref{eq:PFEqua}). 
We revisit a classical result in the following.
\begin{lem} {\textnormal{\textbf{~\cite{wu1982steady}}}}
\label{lem:Existence}
Suppose that the following conditions hold
\begin{enumerate}[(a)]
\item{The microgrid is connected,}
\item{The admittance matrix \textnormal{$\textbf{Y}$} is symmetric,}
\item{$2E_j > E_k$ for all $j,k \in \mathcal{V}$,}
\item{$I_{jk}^{\max} \leq \frac{\pi}{2}B_{jk}$ for all $\{j,k\} \in \mathcal{E}$,}
\item{$\underline{E}_j \big(-B_{jj} + \sum_{k \in \mathcal{V}_I}B_{jk} \big) \geq \big(\sum_{k \in \mathcal{V} \setminus \{j\}}B_{jk}\bar{E}_k\big)$ for all $j \in \mathcal{V}_L$ and the strict inequality holds for at least one $j \in \mathcal{V}_L$,} 
\item{$
P_i \in [\underline{P}_i \; \bar{P}_i] \; \forall i \in \mathcal{V}, Q_i \in [\underline{Q}_i \; \bar{Q}_i] \; \forall i \in \mathcal{V}_L
$,}
\end{enumerate}
where $B_{jk} = Y_{jk}\sin(\phi_{jk})$ and $\underline{P}_i, \bar{P}_i, \;\underline{Q}_i, \bar{Q}_i$ are constants determined by microgrid parameters including line impedance and bounds of voltage regulation. Then there exists a solution to Eq. (\ref{eq:PFEqua}) such that $x \in \mathcal{X}_c$.
\end{lem}
Lemma \ref{lem:Existence} is in fact a direct consequence of Theorem 4 in~\cite{wu1982steady}. Readers are referred to~\cite{wu1982steady} for the proof and details of finding $\underline{P}_i, \bar{P}_i, \;\underline{Q}_i$ and $\bar{Q}_i$. 
\begin{rem}
\normalfont{For simplicity, several conditions in~\cite{wu1982steady} related to loads serviceability are not included in Lemma~\ref{lem:Existence}. Since a transition to island mode is enabled only when the DERs can provide sufficient power to loads in the microgrid, the serviceability requirement is satisfied intrinsically for this work.}
\end{rem}
\begin{rem}
\label{rem:PPSharing}
\normalfont{If the nominal active power injection at inverter buses satisfying $P_i^* \in [\underline{P}_i \; \bar{P}_i], \; \forall i \in \mathcal{V}_I$, Lemma~\ref{lem:Existence} implies $\mathcal{X}_e = \mathcal{X}_S \cap \mathcal{X}_c \neq \emptyset$. Since we focus on microgrid control problem with given $P_i^*$ and $Q_i^*$, we thus assume that $P_i^*$ and $Q_i^*$ selected in the tertiary control layer are chosen such that $P_i^* \in [\underline{P}_i \; \bar{P}_i], \; \forall i \in \mathcal{V}_I$. In this way, $\mathcal{X}_e$ is nonempty if conditions in Lemma~\ref{lem:Existence} hold. We can then focus on designing controller to steer the microgrid to~$\mathcal{X}_e$.}
\end{rem}

\subsection{Distributed Controller Design}
We start our controller design from a simple property of a connected graph $\mathcal{G}_c$. Let $L \in \mathbb{R}^{n_I \times n_I}$ be the Laplacian of $\mathcal{G}_c$. The null space of $L$ is span$\{ {\bf 1}_{n_I} \}$ because $\mathcal{G}_c$ is connected. Observing that ${\bf 1}_{n_I}$ has a close relation with condition (\ref{eq:AcPSharing}), we design $\mu_i(\cdot)$ as a simple linear feedback in terms of $S_{\mathcal{N}_i}$ in the following form
\begin{align}
\label{eq:CoDroop}
\left\{
\begin{array}{ll}
& \dot{x}_i(t) = K_i\sum_{j \in \mathcal{V}_I}
L(i,j)S_j(x(t)) \\
& \text{subj. to } x(t) \in \mathcal{X}_c, \; \forall t\geq 0, 
\end{array}
\right.
\end{align}
where $K_i \in \mathbb{R}^{2 \times 2}$ is the local control gain matrix at bus $i$ to be designed. Define $S_I = \{ S_i, \;i \in \mathcal{V}_I \}$, and $S = \{ S_i, \;i \in \mathcal{V}\}$. Let $K= diag\{K_i, i\in \mathcal{V}_I \}$, and $\bar{L} = L\otimes I_2$. The dynamical model of the microgrid under the proposed inverter control (\ref{eq:CoDroop}) becomes
\begin{align}
\label{eq:CoDroop2}
\left\{
\begin{array}{ll}
& \dot{x}_{I}(t) = K\bar{L}S_{I}(x(t))  \\ 
& \text{subj. to } x(t) \in \mathcal{X}_c, \; \forall t\geq 0, 
\end{array}
\right.
\end{align}
\begin{rem}
\normalfont{The proposed control structure (\ref{eq:CoDroop2}) does not depend on the voltage magnitude information $E_{\mathcal{N}_i}$ that is also available at bus $i$. We will show later that such a control structure is already sufficient to ensure convergence to $\mathcal{X}_e$. In principle, the magnitude information $E_{\mathcal{N}_i}$ can be used to further improve the control performance, especially for voltage regulation. However, we will not study such extension in this paper.}
\end{rem}
Define $\mathcal{O} =\text{span} \{v_p, v_q\}$, where $v_p, v_q \in \mathbb{R}^{2n_I}$, $v_p = [1, 0, 1, \cdots, 0]^T$ and $v_q = [0, 1, 0, \cdots, 1]^T$. The following proposition shows that under some mild conditions, every equilibrium point of system (\ref{eq:CoDroop2}) satisfies the control objectives (\ref{eq:AcPSharing})-(\ref{eq:FReg}).
\begin{prop}
\label{prop:EqStab&PShar}
If Assumption~\ref{assum:Bddxldot} holds and null($K$)$\subseteq \mathcal{O}$, the following statements are equivalent 
\begin{enumerate}[(a)]
\item{The microgrid with dynamics (\ref{eq:CoDroop2}) is in steady state where $\dot{x}_I=0$.}
\item{The desired conditions (\ref{eq:AcPSharing})-(\ref{eq:FReg}) hold.}
\end{enumerate}
\end{prop}
\begin{pf*}{Proof.}
Since null$(L) = {\bf{1}}_{n_I}$, null($\bar{L}$)$= \mathcal{O}$. Null$(K\bar{L}) = \mathcal{O}$ follows directly from null($K$)$\subseteq \mathcal{O}$ and null($\bar{L}$)$= \mathcal{O}$
With null($K\bar{L}$)$= \mathcal{O}$, we have
\begin{align}
\label{eq:EquivalencyStab}
& \dot{x}_I = 0 \text{ such that } x \in \mathcal{X}_c \\ \nonumber \Leftrightarrow & x\in \{x | S_I(x) \in \mathcal{O} \} \cap \mathcal{X}_c  \Leftrightarrow x \in \mathcal{X}_e,
\end{align}
The equivalence between statements (a) and (b) follows from Eq. (\ref{eq:EquivalencyStab}). 
\hfill $\Box$
\end{pf*}
Proposition \ref{prop:EqStab&PShar} reduces the microgrid control problem with numerous requirements to the study of exponential convergence to $\mathcal{X}_e$ of system (\ref{eq:CoDroop2}). We will therefore focus on analyzing system (\ref{eq:CoDroop2}).

\subsection{Analysis of System~(\ref{eq:CoDroop2})}
\label{sec: StabilityAnalysis}
In this subsection, we derive the conditions of exponential convergence to $\mathcal{X}_e$ where all the desired conditions (\ref{eq:AcPSharing})-(\ref{eq:FReg}) follow.  The exponential convergence to $\mathcal{X}_{e}$ of system (\ref{eq:CoDroop2}) is challenging in general due to the nonlinearity of the underlying system and the uncertainty of the load bus states $x_L$. Instead of directly analyzing system (\ref{eq:CoDroop2}), we apply the chain rule to obtain the dynamics of $S_I$ under the proposed control strategy
\begin{align}
\label{eq:DroopInS}
\left\{
\begin{array}{ll}
& \dot{S}_{I}(x(t)) = J_{I,x(t)}K\bar{L}S_{I}(x(t)) + J_{L,x(t)}\dot{x}_{L}(t) \\ 
& x(t) \in \mathcal{X}_c, \; \forall t\geq 0,
\end{array}
\right.
\end{align}
where $J_{I,x}$ and $J_{L,x}$ are the Jacobian matrices of $S_{I}(\cdot)$ evaluated at $x$ with respect to $x_I$ and $x_L$, respectively. Notice that at every time instant and for any $\dot{x}_L$, Eq. (\ref{eq:DroopInS}) describes the dynamics of $S_I$ when the dynamics of $x_I$ is given by Eq. (\ref{eq:CoDroop2}). According to Proposition \ref{prop:EqStab&PShar}, the convergence of $S_I(x(t))$ to $\mathcal{O}$ of system~(\ref{eq:DroopInS}) implies the convergence of the state trajectory $x(t)$ to $\mathcal{X}_e$. The close relation between these two stability properties motivates us to focus on system (\ref{eq:DroopInS}). To simplify notation, we define $z(t) = S_I(x(t))$, $B(t) = J_{I,x(t)}$, and $\mathcal{J}= \{J_{I,x} | x \in \mathcal{X}_c\}$. System (\ref{eq:DroopInS}) can then be written as a linear time varying (LTV) system
\begin{align}
\label{eq:DroopInY}
\left\{
\begin{array}{ll}
& \dot{z}(t) = B(t)K\bar{L}z(t) + w(t)\\ 
& B(t) \in \mathcal{J},
\end{array}
\right.
\end{align}
where $w(t) = J_{L,x(t)}\dot{x}_{L}(t)$ is considered as a disturbance of system (\ref{eq:DroopInY}). With this notation, system (\ref{eq:DroopInY}) becomes a stand alone dynamic system with state variable $z$ subject to unknown disturbance $w(t)$. Note that in system (\ref{eq:DroopInY}), ${w}$ is quadratically bounded by $d(z,\mathcal{O})$
\begin{align*}
||{w}||_2 &= ||J_{L,x}\dot{x}_L||_2 \leq \kappa ||J_{L,x}||_2
 ||\dot{x}_I||_2 \\ 
&= \kappa ||J_{L,x}||_2 ||K\bar{L}z||_2 \leq \zeta ||\bar{L}{z}||_2 \\
&= \zeta\cdot d(z, \mathcal{O}),
\end{align*}
where $\zeta \in \mathbb{R}_+$ is a constant depending on system parameters as well as control gain $K$. Robust stability of the equilibriums of systems with bounded noise was studied in~\cite{vsiljak2000robust}, which is reviewed in the following
\begin{defn}
The set $\mathcal{O}$ is robustly stable of system (\ref{eq:DroopInY}) with degree $\zeta$ if $\mathcal{O}$ is globally exponentially stable for all $w$ such that $||{w}||_2 \leq \zeta d(z, \mathcal{O})$.
\end{defn}
To analyze robust stability for $\mathcal{O}$ of system (\ref{eq:DroopInY}), we apply a standard change of coordinates. Define a change of coordinate matrix 
equation 
\begin{align}
\label{eq:ChangeCor}
T=[v_1, v_2, .., v_{2n-2}, \frac{v_p}{\|v_p\|_2} ,\frac{v_q}{||v_q||_2}], 
\end{align}
where the first $2n-2$ vectors are arbitrary vectors such that $T$ is an orthogonal matrix. Let $\bar{z} = T^{-1}z$ be the state vector in the new coordinate system. The LTV system (\ref{eq:DroopInY}) becomes
\begin{align}
\label{eq:PreReOrderSys}
\left\{
\begin{array}{ll}
& \dot{\bar{z}}(t) = T^{-1}B(t)K\bar{L}T\bar{z}(t) + T^{-1}w(t) \\
& B(t) \in \mathcal{J}.
\end{array}
\right.
\end{align}
Since the last two coordinates of the new basis $V$ span $\mathcal{O}$, the last two column vectors of $\bar{L}T$ are zeros and
\begin{align}
\label{eq:DefnAhat}
T^{-1}B(t)K\bar{L}T = \begin{bmatrix}
\hat{A}_{11}(t) & 0 \\ \hat{A}_{21}(t) & 0
\end{bmatrix},
\end{align}
where $\hat{A}_{11}(\cdot) \in \mathbb{R}^{2(n_I-1) \times 2(n_I-1)}$ and $\hat{A}_{21}(\cdot) \in \mathbb{R}^{2 \times 2(n_I-1)}$. Considering that the dynamics of $d(z(\cdot), \mathcal{O})$ of system (\ref{eq:DroopInY}) is irrelevant to the last two coordinates of the state $\bar{z}$ of system (\ref{eq:PreReOrderSys}), we focus on a reduced order system of (\ref{eq:PreReOrderSys}) with state vector $\hat{z} = [I \; 0]\bar{z} \in \mathbb{R}^{2(n_I-1)}$. Define $\hat{w} = [I \; 0]T^{-1}w \in \mathbb{R}^{2(n_I-1)}$ and $G(B) = [I \; 0]T^{-1}BK\bar{L}T[I \; 0]^T \in \mathbb{R}^{2(n_I-1) \times 2(n_I-1)}$, we have a reduced order system of (\ref{eq:PreReOrderSys})
\begin{align}
\label{eq:ReducedSys}
\left\{
\begin{array}{ll}
& \dot{\hat{z}}(t) = \hat{A}_{11}(t)\hat{z}(t) + \hat{w}(t) \\
& \hat{A}_{11}(t) \in \mathcal{A},
\end{array}
\right.
\end{align}
where $\mathcal{A}:= \{G(B) |\; B \in \mathcal{J}\}$. Similar to $w$ in system (\ref{eq:DroopInY}), $\hat{w}$ is quadratically bounded by the state $\hat{z}$ shown in the following 
\begin{align*}
||\hat{w}||_2 &\leq ||w||_2\leq \zeta ||\hat{z}||_2.
\end{align*}
The following lemma shows that exponential convergence to $\mathcal{X}_e$ of system (\ref{eq:CoDroop2}) follows if the origin of system (\ref{eq:ReducedSys}) is robustly stable. 
\begin{lem}
\label{lem:EqStab}
If the origin of system (\ref{eq:ReducedSys}) is robustly stable with degree $\zeta$, then there exists an non-empty $\mathcal{X}_{c,s} \subseteq \mathcal{X}_c$ such that for all $x(0) \in \mathcal{X}_{c,s}$, $x(t)$ exponentially converges to $\mathcal{X}_e$ for system (\ref{eq:CoDroop2}).
\end{lem} 
\begin{pf*}{Proof.}
Since systems (\ref{eq:ReducedSys}) and (\ref{eq:PreReOrderSys}) share the same dynamics in the space $\mathbb{R}^{2n_I}\setminus \mathcal{O}$, robust stability of the origin of system (\ref{eq:ReducedSys}) implies $\mathcal{O}$ is robustly stable of system (\ref{eq:PreReOrderSys}) with degree $\zeta$. In addition, robust stability of $\mathcal{O}$ of system (\ref{eq:DroopInY}) (or system (\ref{eq:PreReOrderSys}) ) guarantees the trajectory of $z(t)$ is bounded. Define $\mathbb{S}_I = \{S_I(x)|x \in \mathcal{X}_c \}$. With the bounded trajectory of $z(t) = S_I(x(t))$, there exists $\mathcal{X}_{c,s} \subseteq \mathcal{X}_c$ such that for all $x(0) \in \mathcal{X}_{c,s}$, $S_I(x(t)) \in \mathbb{S}_I \; \forall t$, which implies $x(t) \in \mathcal{X}_c \; \forall t$. Therefore, for all initial $x(0) \in \mathcal{X}_{c,s}$, $z(t)=S_I(x(t))$ converges to $\mathcal{O}$ with $x(t) \in \mathcal{X}_c$ for all time in system (\ref{eq:DroopInY}) if system (\ref{eq:ReducedSys}) is robustly stable. Recall that Eq.~(\ref{eq:DroopInY}) describes the dynamics of $S_I(\cdot)$ when the dynamics of $x_I$ is given by Eq.~(\ref{eq:CoDroop2}). We then conclude that for all $x(0) \in \mathcal{X}_{c,s}$, $x(\cdot)$ exponentially converges to $\mathcal{X}_e$ for system~(\ref{eq:CoDroop2}) if the origin of system~(\ref{eq:ReducedSys}) is robustly stable. \hfill $\Box$
\end{pf*}
We now provide a set of sufficient conditions for robust stability of the origin of system (\ref{eq:ReducedSys})
\begin{prop}
\label{prop:RobStab}
The origin of system (\ref{eq:ReducedSys}) is robustly stable with degree $\zeta$ if there exist $\epsilon,\xi \in \mathbb{R}_+$, $U = U^T \succ 0$ such that
\begin{align}
\label{eq:ConForStability}
\begin{bmatrix}
\hat{A}_{11}^TU + U\hat{A}_{11}+ \epsilon\zeta I+\xi U & U \\ U & - \epsilon I
\end{bmatrix}\preceq 0 ,
\end{align}
for all $\hat{A}_{11}\in \mathcal{A}$.
\end{prop}
\begin{pf*}{Proof.}
The proof is similar to linear time invariant system case discussed in~\cite{vsiljak2000robust}. Eq.~(\ref{eq:ConForStability}) can be derived by quadratic Lyapunov function argument. If there exist a Lyapunov function $V(\hat{z}) = \hat{z}^TU\hat{z}$ such that for all $\hat{A}_{11} \in \mathcal{A}$, $\dot{V}(\hat{z}) \leq -\xi{V}(\hat{z}) $, then the origin of system~(\ref{eq:ReducedSys}) is exponentially stable. The conditions for $\dot{V}(\hat{z}) \leq -\xi{V}(\hat{z})$, $\forall \hat{A}_{11} \in \mathcal{A}$ such that $||w||_2\leq \zeta ||\hat{z}||_2$ are shown in the following
\begin{align*}
&\dot{V}(\hat{z}) \leq -\xi V \text{ s.t. }||w||_2\leq \zeta ||\hat{z}||_2 \\
\Longleftrightarrow & \begin{bmatrix}
\hat{z}^T \\ \hat{w}^T
\end{bmatrix}\begin{bmatrix}
\hat{A}_{11}^TU +U\hat{A}_{11} + \xi U & U \\ U & 0 
\end{bmatrix}\begin{bmatrix}
\hat{z} \\ \hat{w}
\end{bmatrix} \leq 0 \\
& \text{ s.t. } \begin{bmatrix}
\hat{z}^T \\ \hat{w}^T
\end{bmatrix} \begin{bmatrix}
-\zeta I & 0 \\ 0 & I 
\end{bmatrix} \begin{bmatrix}
\hat{z} \\ \hat{w}
\end{bmatrix} \leq 0 \\
\Longleftrightarrow & \begin{bmatrix}
\hat{z}^T \\ \hat{w}^T
\end{bmatrix}\begin{bmatrix}
\hat{A}_{11}^TU +U\hat{A}_{11}+\epsilon\zeta I + \xi U & U \\ U & -\epsilon I 
\end{bmatrix}\begin{bmatrix}
\hat{z} \\ \hat{w}
\end{bmatrix} \leq 0.
\end{align*}
for all $\hat{A}_{11}\in \mathcal{A}$. S-procedure is applied for the last step, which completes the proof. \hfill $\Box$
\end{pf*}
Note that if $\hat{A}_{11}$ is polytopic, the condition (\ref{eq:ConForStability}) can be formulated into bilinear matrix inequalities (BMIs). The condition can then be checked numerically. However, $\hat{A}_{11}$ is not polytopic, so we will instead develop a way to find a convex set containing $\hat{A}_{11}$ in the next subsection. 

With Lemma~\ref{lem:EqStab} and Proposition~\ref{prop:RobStab}, conditions of exponential convergence to $\mathcal{X}_{e}$ of (\ref{eq:CoDroop2}) can be obtained:
\begin{thm}
\label{thm: Key}
If Assumption~\ref{assum:Bddxldot}, hypotheses in Lemma \ref{lem:Existence} and Eq. (\ref{eq:ConForStability}) hold, then there exists an non-empty $\mathcal{X}_{c,s} \subseteq \mathcal{X}_c$ such that for $x(0) \in \mathcal{X}_{c,s}$, the microgrid (\ref{eq:CoDroop2}) converges exponentially to the set $\mathcal{X}_{e}$ where the control objectives including (\ref{eq:AcPSharing})-(\ref{eq:FReg}) are all satisfied. 
\end{thm}
\begin{pf*}{Proof.}
The origin of system (\ref{eq:ReducedSys}) is robustly stable due to Proposition~\ref{prop:RobStab}. By Lemma~\ref{lem:EqStab}, robust stability of system (\ref{eq:ReducedSys}) implies the existence of $\mathcal{X}_{c,s} \subseteq \mathcal{X}_c$ such that for all $x(0) \in \mathcal{X}_{c,s}$, the trajectories converge to $\mathcal{X}_e$ of system (\ref{eq:CoDroop2}). Since Eq.~(\ref{eq:ConForStability}) ensures null($K$)$\subseteq \mathcal{O}$, $\mathcal{X}_e$ is equivalent to desired control objectives (\ref{eq:AcPSharing})-(\ref{eq:FReg}) from Proposition~\ref{prop:EqStab&PShar} .  \hfill $\Box$
\end{pf*}
\normalfont{The result of Theorem \ref{thm: Key} is robust with respect to small variations of system parameters. As long as the perturbations of the admittance matrix are small enough such that $J_{I,x} \in \mathcal{J}$, the exponential convergence to $\mathcal{X}_{e}$ for some $x(0) \in \mathcal{X}_{c,s}$ still follows from Theorem \ref{thm: Key}. Different from most of the literature, the proposed controller can be applied to mixed $R/X$ ratio distribution lines and general microgrid topology including acyclic and mesh networks. Furthermore, the controller can meet all the main control objectives without the separation of time scale, which distinguishes it from the mainstream droop control methods.}

\subsection{Feedback Gain Design}
\label{sec:RealWorldImple}
In this subsection, we propose a constructive way to find a feedback gain $K$ satisfying Eq.~(\ref{eq:ConForStability}). The difficulty lies in checking the feasibility Eq.~(\ref{eq:ConForStability}). As discussed in the last subsection, the robust stability condition in Eq. (\ref{eq:ConForStability}) can not be directly formulated into BMIs because $\mathcal{A}$ is not polytopic. We will first derive a convex hull containing $\mathcal{A}$ by analyzing the Jacobian of the power flow equations (\ref{eq:PFEqua}) so that Eq.~(\ref{eq:ConForStability}) can be checked by solving several BMIs. Secondly, instead of only checking the feasibility, we formulate the BMIs into an optimization problem to enhance the robustness. 

Define $C_i \subset \mathcal{E}$ and $\mathcal{X}_{C_i}$ as a polyhedron replacing the constraint $|\theta_i -\theta_j| \leq \gamma$ in $\mathcal{X}_c$ by $\theta_i -\theta_j = 0$ for all $\{i,j\} \in C_i$. Let $\Phi:= \{\phi_{ij} | \; |\phi_{ij}|+\gamma \leq \frac{\pi}{2} \}$. An approximated convex hull can be found by the following Proposition.
\begin{prop}
\label{prop:LMIs}
If every entry of the admittance matrix \normalfont{Y} satisfies $\phi_{ij} \in \{\pm \frac{\pi}{2}, \Phi \}$, the upper and lower bounds of every entry $(i,j)$ of $J_{I,x} \in \mathcal{J}$ are 
\begin{align}
\label{eq:Jaco_bdd}
& \bar{J}_{I,x}(i,j) = \max\{J_{I,z}(i,j), z \in Z\},  \\ \nonumber 
& \underline{J}_{I,x}(i,j) = \min\{J_{I,z}(i,j), z \in Z\},
\end{align}
where $Z = \bigcup_{C_i \in 2^{\mathcal{E}}}v(\mathcal{X}_{C_i})$.
\end{prop}
\begin{pf*}{Proof.}
According to the power flow equation (\ref{eq:PFEqua}), the derivative of the active power injection at bus $i$ with respect to different variants are 
\begin{align}
\label{eq:Jaco_Elem}
\begin{split}
& \frac{\partial P_i(x)}{\partial \theta_i} = -E_i\sum_{j \in \mathcal{V} \setminus i}Y_{ij}E_j\sin(\theta_i -\theta_j -\phi_{ij}), \\ 
& \frac{\partial P_i(x)}{\partial E_i} = 2E_iY_{ii} \cos( -\phi_{ii}) + \sum_{j \in \mathcal{V}\setminus i} Y_{ij}E_j \text{cos}(\theta_i - \theta_j - \phi_{ij}),  \\
& \frac{\partial P_i(x)}{\partial \theta_j} = E_iY_{ij}E_j\sin(\theta_i -\theta_j -\phi_{ij}) , \; j \neq i, \\
& \frac{\partial P_i(x)}{\partial E_j} = E_iY_{ij}\cos(\theta_i -\theta_j -\phi_{ij}), \;  j \neq i.
\end{split}
\end{align}
Since for all entries in the admittance matrix Y, $\phi_{ij} \in \{\pm \frac{\pi}{2}, \Phi \}$, every summand of the Jacobian elements in (\ref{eq:Jaco_Elem}) has the maximal and minimal points at $x \in Z$ if $x \in \mathcal{X}_c$. In addition, every summand corresponds to different set of $E_j$ and $\theta_i - \theta_j$ in each Jacobian element in (\ref{eq:Jaco_Elem}), the Jacobian elements have the maximal and minimal points at $x \in Z$ if $x \in \mathcal{X}_c$. The same conclusion for the reactive part is reached by a similar argument. \hfill $\Box$
\end{pf*}
Given the upper and lower bounds of every entry of $J_{I,x} \in \mathcal{J}$, a convex hull containing $\mathcal{J}$ can be found. Let $\mathbb{F}$ be a set $\{-1,1\}$. Define $\varrho_i \in \mathbb{F}^{2n_I \times 2n_I}$ such that $\cup_{i\in [l]}\varrho_{i} = \mathbb{F}^{2n_I \times 2n_I}$, where $l = 2^{4n_I^2}$. Define $D_{\varrho_i} \in \mathbb{R}^{2n_I \times 2n_I}$ such that
\begin{align*}
& D_{\varrho_i}(j,k) = \bar{J}_{I,x}(j,k) \text{ if } \varrho_i(j,k)= 1, \\ 
& D_{\varrho_i}(j,k) = \underline{J}_{I,x}(j,k) \text{ if } \varrho_i(j,k)= -1.
\end{align*}
The following lemma is a simple consequence of Proposition 3 and provides a convex hull containing $\mathcal{J}$.
\begin{lem}
\label{lem: ConvexHull}
If every entry of the admittance matrix \normalfont{Y} satisfy $\phi_{ij} \in \{\pm \frac{\pi}{2}, \Phi \}$, the convex hull 
\begin{align*}
\mathcal{\bar{D}}:= CO\{D_{\varrho_1}, ..., D_{\varrho_l}\}
\end{align*}
contains $\mathcal{J}$.
\end{lem}
The results of Lemma~\ref{lem: ConvexHull} allows us to replace Eq. (\ref{eq:ConForStability}) in Theorem~\ref{thm: Key} by BMIs. Instead of only finding $K$ such that Eq.~(\ref{eq:ConForStability}) is feasible, we propose to design $K$ by solving the following optimization problem subject to BMI constraints
\begin{equation}
\begin{aligned} 
\label{eq: LMIsStability}
\mathop{\text{maximize}}_{K,U,\epsilon, \xi} \;
 & \xi \\
\text{subject to } & M_i(\epsilon,\zeta,\xi) \preceq 0,
\end{aligned}
\end{equation}
where $\epsilon, \xi \in \mathbb{R}_+$ and
\begin{align*}
& \hat{A}_{11,i} = [I, 0]T^{-1}D_{\varrho_i}K\bar{L}T[I, 0]^T, \\
& M_i(\cdot) = \begin{bmatrix}
\hat{A}_{11,i}^TU + U\hat{A}_{11,i} + \epsilon\zeta I + \xi U & U \\ U & - \epsilon I
\end{bmatrix}, \\
& U = U^T \succ 0,\; i \in [l].
\end{align*}
\begin{rem}
\label{rem:K_bounded}
\normalfont{The maximization of $\xi$ is for the purpose of improving the convergence rate.  Notice that maximizing $\xi$ is only meaningful if an upper bound of $||K||_2$ is imposed. More discussions on the upper bound will be included in the next section.}
\end{rem}
The following corollary is a direct consequence of Eq.~(\ref{eq: LMIsStability}) and Theorem~\ref{thm: Key}
\begin{cor}
\label{cor: KeyLMI}
If Assumption~\ref{assum:Bddxldot}, the hypotheses in Lemma \ref{lem:Existence} hold, and the optimal solution $K$ for Eq.~(\ref{eq: LMIsStability}) is implemented in~Eq.~(\ref{eq:CoDroop2}), then there exists an non-empty $\mathcal{X}_{c,s} \subseteq \mathcal{X}_c$ such that for $x(0) \in \mathcal{X}_{c,s}$, the microgrid (\ref{eq:CoDroop2}) converges exponentially to the set $\mathcal{X}_{e}$ where the control objectives including (\ref{eq:AcPSharing})-(\ref{eq:FReg}) are all satisfied. 
\end{cor}
We want to comment that the convex optimization problem~(\ref{eq: LMIsStability}) subject to BMI constraints is NP-hard to solve in general. However, efficient algorithms~\cite{hassibi1999path} and \cite{kanev2004robust} are available if an initial feasible solution can be found. We can first substitute $U$ by some simple positive definite matrices as an initial guess and check the feasibility. If $U$ is a feasible solution of Eq.~(\ref{eq: LMIsStability}), algorithms~\cite{hassibi1999path} or \cite{kanev2004robust} can be applied to find a local optimal solution. Those algorithms only involve several linear matrix inequalities (LMIs) instead of BMIs, which can be effectively solved by using various existing convex optimization algorithms~\cite{cvx,gb08}.
\begin{rem}(Trade-off between Complexity and Robustness)
\label{rem:ComplexRobust}
\normalfont{Let $\bar{\mathcal{J}}:= CO\{J_{I,z}, z\in Z\}$. The convex hull $\bar{\mathcal{J}}$ does not necessarily contain $\mathcal{J}$. However, there exists a compact set $\mathcal{X}_{c,s} \subseteq \mathcal{X}_{c}$ such that $\bar{\mathcal{J}} \supseteq \{ J_{I,x}, x \in \mathcal{X}_{c,s} \}$. Hence, if one solves Eq.~(\ref{eq: LMIsStability}) by replacing $v(\bar{\mathcal{D}})$ with $v(\bar{\mathcal{J}})$, Corollary~\ref{cor: KeyLMI} still implies exponential convergence to a subset of $\mathcal{X}_{e}$ for some initial $x(0)$. The benefit of replacing $v(\bar{\mathcal{D}})$ with $v(\bar{\mathcal{J}})$ is reduction of complexity, where $|v(\bar{\mathcal{D}})| \sim 2^{4n_I^2}$ while $|v(\bar{\mathcal{J}})| \sim 2^{n}\cdot3^{|\mathcal{E}|}$. The difference becomes evident for large $n_I$ and one may prefer to solve the optimization problem~(\ref{eq: LMIsStability}) by substituting $v(\bar{\mathcal{D}})$ to $v(\bar{\mathcal{J}})$ for microgrids with larger $n_I$.}
\end{rem}

\begin{figure}
\begin{center}
\centerline{\scalebox{0.85}{\includegraphics{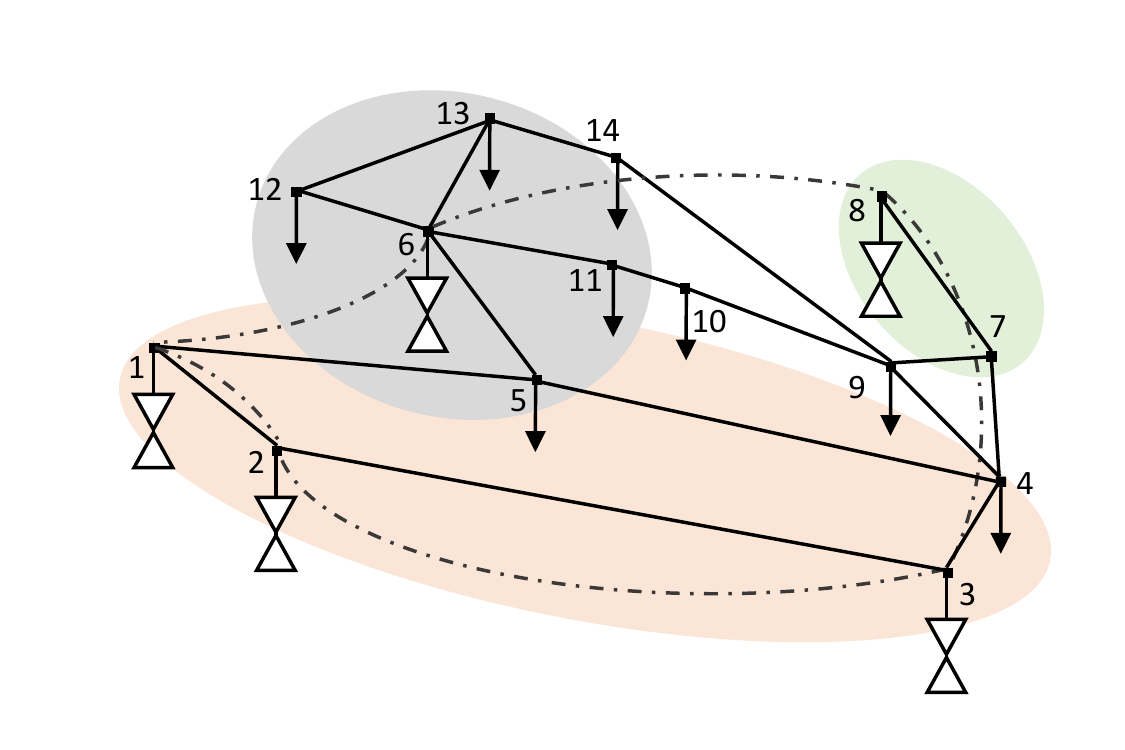}}}
\caption{IEEE 14 bus system. The dash lines are the communication links. Every generation bus, possibly with local load, is replaced by one VSI in this work for simplicity.
}
\label{fig:14busTopo}
\end{center}
\end{figure}

As discussed in Remark~\ref{rem:ComplexRobust}, the number of the BMI constraints in the optimization problem~(\ref{eq: LMIsStability}) grows exponentially with respect to $n_I$. Finding $K$ for microgrids with large number of bus can become challenging even one replacing $v(\bar{\mathcal{D}})$ by $v(\bar{\mathcal{J}})$. Here, we observe that matrices in ${\mathcal{J}}$ are block diagonal. The number of blocks equals to the number of groups of inverter bus separated by load bus. For example, all matrices of ${\mathcal{J}}$ of IEEE 14 bus system shown in Fig.~\ref{fig:14busTopo} has three diagonal blocks. One of the blocks corresponds to $1^{\text{th}}-3^{\text{rd}}$ inverter buses, and the values of the entries depend on the states of $1^{\text{th}}-5^{\text{th}}$ buses, which are covered with the same color in Fig.~\ref{fig:14busTopo}. Similar argument applies to the other two blocks. 

Following the motivating example above, we introduce several notations to define properties of every block in $\mathcal{J}$. Let $c_i$ be the set of inverter buses associated with block $i$. The dimension of each block is then given by $|c_i|$. Similar arguments to Proposition~3 and Lemma~3 are made to define the convex hull $\bar{\mathcal{D}}_i$ associated with block $i$. The number of vertices of $\bar{\mathcal{D}}_i$ is then $l_i = 2^{4|c_i|^2}$. Notice that the large number of BMIs in Eq.~(\ref{eq: LMIsStability}) is originated from the number of vertices of the convex hull $\bar{\mathcal{D}}$. In the following lemma, we will show that each block in the block diagonal Jacobian matrix can be viewed separately in finding a feasible solution of Eq.~(\ref{eq: LMIsStability}), resulting in much less number of BMIs involved.

\begin{lem}
\label{lem:PolynomialBMI}
Every $K$ that satisfies Eq.~(\ref{eq:PolynomialBMI}) is a feasible solution of Eq.~(\ref{eq: LMIsStability}).
\begin{align}
\label{eq:PolynomialBMI}
&\underline{\lambda}(D_{\varrho_k,i}K_{c_i}+K_{c_i}^TD_{\varrho_k,i}^T) \leq -d < 0 \\ \nonumber
& \forall k \in [l_i], \;\;
\forall i \in [n_c],
\end{align}
where $K_{c_i} = diag\{K_j, j\in c_i \}$ and $n_c$ is the number of diagonal block for matrices in $\mathcal{J}$.
\end{lem}
\begin{pf*}{Proof.}
From Eq.~(\ref{eq:PolynomialBMI}), there exist $U = U^T\succ 0$ such that 
\begin{align*}
\underline{\lambda}(\hat{A}_{11,i}^TU+U\hat{A}_{11,i}) \leq -\bar{d} < 0, \; \forall i \in [l].
\end{align*}
By Schur compliment, $M_i(\epsilon,\zeta,\xi) \preceq 0 $ if and only if
\begin{align}
\label{eq:Proof_feas}
(\hat{A}_{11,i}^TU + U\hat{A}_{11,i} + \xi U)+ \epsilon(UU + \zeta I) \preceq 0. 
\end{align}
We can find $\epsilon,\xi \in \mathbb{R}_+$ such that Eq.~(\ref{eq:Proof_feas}) holds when $\bar{d} \in \mathbb{R}_+$ and the desired result follows.  \hfill $\Box$
\end{pf*}
Lemma~\ref{lem:PolynomialBMI} allows us to find $K$ by solving Eq.~(\ref{eq:PolynomialBMI}) instead of Eq.~(\ref{eq: LMIsStability}). The computational complexity is greatly reduced because the number of BMIs associated to Eq.~(\ref{eq:PolynomialBMI}) increases with respect to $\bar{l}=\max_{i \in [n_c]}l_i$ instead of $n$. In most microgrid networks, $\bar{l}$ is much smaller than $n$.
Solving Eq.~(\ref{eq:PolynomialBMI}) instead of Eq.~(\ref{eq: LMIsStability}) enhances the practicability of the proposed distributed control method, especially for large-scale microgrids.

\section{Flexible Operation of Microgrids}
\label{sec:FlexOperMicro}
In this section, we study the plug-and-play feature of DERs under the proposed distributed control framework. The robustness under communication failure is also analyzed to solidify the proposed distributed controller. 

\subsection{Plug-and-Play Feature}
The plug-and-play feature of the DERs refers to the property that one DER can be plugged or unplugged to a microgrid without re-engineering the entire control. We consider a general case where part of inverters may be disconnected from the microgrid abruptly due to some severe events. Let $\mathcal{V}_{I_f} \subset \mathcal{V}_{I}$ be a set of normal operating inverter buses so that $\mathcal{V}_I \setminus \mathcal{V}_{I_f}$ is the set of disconnected inverters. The voltage magnitude and phase angle dynamics at buses $ i \in \mathcal{V}_I \setminus \mathcal{V}_{I_f}$ become unknown and are categorized as load buses. For this reason, we will partition the buses by $\mathcal{V} = \mathcal{V}_{I_f} \cup \mathcal{V}_{L_f}$, for which $\mathcal{V}_{I_f}$ contains all the normal operating inverters, while $\mathcal{V}_{L_f} \triangleq \mathcal{V} \setminus \mathcal{V}_{I_f}$  consists of all the other buses including load buses $V_L$ or disconnected inverter buses $\mathcal{V}_I\setminus \mathcal{V}_{I_f}$. 
Consider the case that the communication network remains intact when some inverters are disconnected from the microgrid. The communication graph between the operating inverters is defined as $\mathcal{G}_{f} := (\mathcal{V}_{I_f}, \mathcal{E}_{f})$, where $\mathcal{E}_{f} = \mathcal{E}_{c} \setminus \{ \{i,j\} | i \in \mathcal{V}_I \setminus \mathcal{V}_{I_f} \}$. Let $L_{f}$ be the Laplacian of $\mathcal{G}_{f}$. The control law of the ``fault'' microgrid is reduced from Eq.~(\ref{eq:CoDroop2}) to
\begin{align}
\label{eq:CoDroopLDG}
\left\{
\begin{array}{ll}
& \dot{x}_i(t) = K_i\sum_{j \in \mathcal{V}_{I_f}}
L_{f}(i,j)S_i(x(t)) \\
& \text{subj. to } x(t) \in \mathcal{X}_c, \; \forall t\geq 0.
\end{array}
\right.
\end{align}
Notice that due to the assumption on intact communication network, the control law can be autonomously transformed to Eq. (\ref{eq:CoDroopLDG}) in response to the change of microgrid operating conditions. Let $x_{I_f}:= \{x_i, i \in \mathcal{V}_{I_f} \}$ and $S_{I_f}:= \{S_i, i \in \mathcal{V}_{I_f} \}$. The microgrid dynamics with the controller (\ref{eq:CoDroopLDG}) is rewritten as follows
\begin{align}
\label{eq:LDGsDroop}
\left\{
\begin{array}{ll}
& \dot{x}_{I_f}(t) = K_{I_f}\bar{L}_{f}S_{I_f}(x(t)) \\
& \text{subj. to }x(t) \in \mathcal{X}_c, \; \forall t\geq 0,
\end{array}
\right.
\end{align}
where $K_{I_f}= diag\{K_i, i\in \mathcal{V}_{I_f} \}$, and $\bar{L}_f = L_f \otimes I_2$. Denoted $\mathcal{O}_f$ as the proportional power sharing space of the inverters $i \in \mathcal{V}_{I_f}$. If $\mathcal{G}_{f}$ stays connected, the null space of $\bar{L}_{f}$ remains equivalent to $\mathcal{O}_f$. The ``reduced'' microgrid (\ref{eq:LDGsDroop}) can therefore be analyzed through a similar way discussed in the last section. Define $\bar{\mathcal{D}}_f = CO(D_{\varrho_{1_f}}, \cdots, D_{\varrho_{l_f}})$ as Jacobians of $S_{I_f}$ when $x \in \mathcal{X}_c$ and $l_f$ as the number of vertices of $\bar{\mathcal{D}}_f$. The control objectives (\ref{eq:AcPSharing})-(\ref{eq:FReg}) of the reduced microgrid follow if all the hypotheses in Corollary \ref{cor: KeyLMI} hold except that Eq.~(\ref{eq: LMIsStability}) is replaced by the following 
\begin{equation}
\label{eq: LMIsStabilityFailDG}
\begin{aligned} 
\mathop{\text{maximize}}_{\epsilon_f,K_f,U_f,\xi_f}\; & \xi_f \\ 
\text{subject to } & M_{i,f}(\epsilon_f,\zeta_f,\xi_f) \preceq 0,
\end{aligned}
\end{equation}
where $\epsilon_f, \xi_f \in \mathbb{R}_+$ and
\begin{align*}
& \hat{A}_{11,i,f} = [I, 0]T_f^{-1}D_{\varrho_{i_f}}K_{I_f}\bar{L}_{f}T_f[I, 0]^T, \\
& M_{i,f}(\cdot) = \begin{bmatrix}
\hat{A}_{11,i,f}^TU_f + U_f\hat{A}_{11,i,f} +\epsilon_f\zeta_f I + \xi_f U_f&  U_f \\ U_f & -\epsilon_f I
\end{bmatrix}, \\
& U_f = U_f^T \succ 0,\; i \in [l_f].
\end{align*}
The orthogonal matrix $T_f$ of system (\ref{eq:LDGsDroop}) is defined similar to $T$ as shown in Eq.~(\ref{eq:ChangeCor}), where the last two coordinates span the space of proportional power sharing. If $K$ is found such that BMI constraints of Eq.~(\ref{eq: LMIsStability}) and (\ref{eq: LMIsStabilityFailDG}) are both satisfied for $\mathcal{V}_{I}$ and all possible $\mathcal{V}_{I_f}$ respectively, the proposed controller has the plug-and-play feature. 

In fact, we will show a lemma proving that the Eq.~(\ref{eq: LMIsStabilityFailDG}) is feasible if Eq.~(\ref{eq: LMIsStability}) is feasible and $\mathcal{G}_{f}$ is connected. 
\begin{lem}
\label{lem:SubmatrixStability}
If $\mathcal{G}_{f}$ is connected and $\epsilon^*,K^*, U^*, \xi$ are feasible solution of Eq.~(\ref{eq: LMIsStability}), then by setting $K_{I_f}=diag\{K_i^*, i \in \mathcal{V}_{I_f} \}$, BMI constraints in Eq.~(\ref{eq: LMIsStabilityFailDG}) reduce to LMIs and are feasible.
\end{lem}
\begin{pf*}{Proof.}
By definition of the Laplacian matrix, $\bar{L}$ and $\bar{L}_f$ are positive semidefinite. In addition, because $\mathcal{G}_c$ and $\mathcal{G}_f$ are connected, 
\begin{align}
\label{eq:OrgLaplacianPos}
& T^{-1}\bar{L}T = \begin{bmatrix}
\bar{L}_1 & 0 \\ 0 & 0
\end{bmatrix}, \;\; 0 \prec \bar{L}_1 \in \mathbb{R}^{2(n_I-1) \times 2(n_I-1)}, \\
\label{eq:LaplacianPos}
& T^{-1}_f\bar{L}_fT_f = \begin{bmatrix}
\bar{L}_{f1} & 0 \\ 0 & 0
\end{bmatrix}, \;\; 0 \prec \bar{L}_{f1} \in \mathbb{R}^{2(|\mathcal{V}_{I_f}|-1) \times 2(|\mathcal{V}_{I_f}|-1)}.
\end{align}
Since $\epsilon^*, K^*, \;U^*$ solve Eq.~(\ref{eq: LMIsStability}), the following condition holds for all $i \in [l]$
\begin{align}
\label{eq:MInEq}
&[U^*, 0]T^{-1}D_{\varrho_i}K^*T[\bar{L}_1, 0]^T + \\ \nonumber
& \hspace{2mm}[\bar{L}_1, 0]T^{-1}(K^*)^TD^T_{\varrho_i}T[U^*, 0]^T
\prec -(\epsilon U^*U^*+\xi^* U^*).
\end{align}
The matrix inequalities are derived by applying Schur complement to BMI constraints of Eq.~(\ref{eq: LMIsStability}) and symmetric property of $\bar{L}_1$. Due to Eq.~(\ref{eq:MInEq}), the eigenvalues of $H_i:= D_{\varrho_i}K^* + (K^*)^TD_{\varrho_i}^T$ can be transformed into the following form by spectral decomposition
\begin{align*}
\bar{H}_i = \begin{bmatrix}
\bar{H}_{i1} & 0 \\ 0 & \bar{H}_{i2}
\end{bmatrix}, 
\end{align*}
where $\lambda(\bar{H}_{i1})\subset \mathbb{R}_{-c}, c>0$. The matrix $\bar{H}_{i1}$ is a bijective mapping from $\mathcal{H}_i$ to $\mathcal{H}_i$, where $\mathcal{H}_i$ is some Euclidean space such that $\mathcal{H}_i \supseteq \mathbb{R}^{2n_I} \setminus \mathcal{O}$. Notice that $D_{\varrho_{i_f}}K_{I_f} +K_{I_f}^TD_{\varrho_{i_f}}^T$ is a principal submatrix of $\bar{H}_i$ and
\begin{align*}
H_{r,i}:= [I, 0]T^{-1}_f(D_{\varrho_i,f}K_{I_f} +K_{I_f}^TD_{\varrho_i,f}^T )T_f[I, 0]^T
\end{align*}
is a linear mapping from $\mathbb{R}^{2|\mathcal{V}_{f_I}|} \setminus \mathcal{O}_f$ to $\mathbb{R}^{2|\mathcal{V}_{f_I}|} \setminus \mathcal{O}_f$, where $\mathbb{R}^{2|\mathcal{V}_{f_I}|} \setminus \mathcal{O}_f \subset \mathcal{H}_i$. As a result, $H_{r,i}$ is a principal submatrix of $\bar{H}_{i1}$ with a proper basis of $\mathcal{H}_i$. By the Cauchy interlace theorem, ${\lambda}(H_{r,i}) \subset \mathbb{R}_{-c}$ for all $i \in [l_f]$. The condition ${\lambda}(H_{r,i}) \subset \mathbb{R}_{-c}$ is sufficient for the feasibility of BMI constraints in Eq.~(\ref{eq: LMIsStabilityFailDG}) by Schur compliment and the property of $\bar{L}_{f}$ shown in Eq. (\ref{eq:LaplacianPos}).
  \hfill $\Box$
\end{pf*}
With Lemma~\ref{lem:SubmatrixStability}, we can conclude stability of the ``fault'' microgrids by only solving optimization problem~(\ref{eq: LMIsStability}), which reduces the computational efforts to find a proper $K$ for various microgrids operating conditions.

\subsection{Loss of Communication links}
We consider the case of communication failures in this subsection. Let $\mathcal{G}_{l} := (\mathcal{V}_I, \mathcal{E}_{l})$ be a simple communication graph when communication links $\mathcal{E}_c\setminus \mathcal{E}_{l}$ are failed. Denote $L_l$ as the Laplacian of $\mathcal{G}_{l}$. We can write down the microgrid dynamics with the controller (\ref{eq:CoDroop2}) by replacing $\bar{L}$ by $\bar{L}_l$. The analyses in the Section \ref{sec:ContrDesign} apply and we can derive conditions to ensure objectives (\ref{eq:AcPSharing})-(\ref{eq:FReg}) under loss of communication links.

Consider a more general scenario when both part of the DERs and communication links failed. Let $\mathcal{G}_{fl} = (\mathcal{V}_{I_f}, \mathcal{E}_{fl})$ and $L_{fl}$ be the graph and Laplacian under the operating condition, respectively. Assume the communication links
\begin{align}
\label{eq:LossDGCommReq}
\{ (i,j), i \in \mathcal{V}_{I_f}, j \in \mathcal{V}_I \setminus \mathcal{V}_{I_f} \} \subseteq \mathcal{E} \setminus \mathcal{E}_{fl}
\end{align}
remains intact so that the inverters can still update the Laplacian $L_{fl}$ autonomously. The desired properties in the scenario of loss of DERs and communication links tie closely to the feasibility of Eq.~(\ref{eq: LMIsStabilityFailDG}). We conclude this section by the following theorem
\begin{thm}
\label{thm: Key2}
If the hypotheses in Corollary \ref{cor: KeyLMI} hold, then there exists an non-empty $\mathcal{X}_{c,s} \subseteq \mathcal{X}_c$ such that for some given loads and initial $x(0) \in \mathcal{X}_{c,s}$, the microgrid converges exponentially to the set where the control objectives including (\ref{eq:AcPSharing})-(\ref{eq:FReg}) are all satisfied if any one of the following operating conditions holds. 
\begin{enumerate}[(a)]
\item{Microgrid operates normally with the dynamics (\ref{eq:CoDroop2}),}
\item{The DERs $i \in \mathcal{V}_I \setminus \mathcal{V}_{I_f}$ are disconnected with the microgrid, where $\mathcal{G}_{f}$ remains connected. The microgrid dynamics is described by (\ref{eq:LDGsDroop}). Voltage dynamics of buses $\mathcal{V} \setminus \mathcal{V}_{I_f}$ is bounded in the form of Proposition~\ref{assum:Bddxldot},}
\item{The communication network is reduced to a connected simple graph $\mathcal{G}_{l}$. The microgrid dynamics is described by Eq. (\ref{eq:CoDroop2}) with $\bar{L}$ replaced by $\bar{L}_l$,}
\item{Both parts of the DERs are unplugged and communication links are failed. The reduced communication graph $\mathcal{G}_{fl}$ stays connected. The dynamics is described by Eq. (\ref{eq:LDGsDroop}) with $\bar{L}_{f}$ replaced by $\bar{L}_{fl}$. Voltage dynamics of buses $\mathcal{V} \setminus \mathcal{V}_{I_{fl}}$ is bounded in the form of Proposition~\ref{assum:Bddxldot}.}
\end{enumerate}
\end{thm}
\begin{pf*}{Proof.}
Every microgrid condition listed above has its own reduced order system in the form of Eq. (\ref{eq:ReducedSys}). To prove the theorem, it is sufficient to show that all the reduced order systems are robustly stable with the hypotheses of Corollary~\ref{cor: KeyLMI}. Namely, the robust stability condition of every reduced order system in the form of Eq.~(\ref{eq: LMIsStability}) is feasible. The operating conditions (a) and (b) have their robust stability conditions hold due to Corollary~\ref{cor: KeyLMI} and Lemma~\ref{lem:SubmatrixStability}. The operating condition~(c) has the same dynamics to Eq. (\ref{eq:CoDroop2}), only $\bar{L}$ is replaced by $\bar{L}_l$. Notice that the eigenspace of $0$ for both $\bar{L}$ and $\bar{L}_l$ are the same. As a result, when $\bar{L}$ is replaced by $\bar{L}_l$ in Eq.~(\ref{eq: LMIsStability}), the new BMI constraints remain feasible with the gain $K$ optimizing Eq.~(\ref{eq: LMIsStability}). The robust stability conditions of operating scenario~(c) follow. The feasibility of BMI constraints for operating scenario~(d) follows if Eq.~(\ref{eq: LMIsStabilityFailDG}) is feasible because the only difference between the two is the replacement of $\bar{L}_f$ with $\bar{L}_{fl}$, where similar arguments from scenario (a) to (c) apply. \hfill $\Box$
\end{pf*}  
\begin{rem} 
\normalfont{All the operating conditions listed in Theorem \ref{thm: Key2} require that the communication graph stays connected. Otherwise, the null space of the Laplacian matrix is bigger than the space of proportional power sharing. Proportional power sharing is only guaranteed between inverters which are connected by the remaining communication network. Toward this end, tertiary control may help establish proportional power sharing. The normalized power difference between ``disconnected'' components is first collected. The tertiary controller then broadcasts a new voltage reference to offset the power difference. Proportional power sharing can be gradually established with an aid of tertiary control.}
\end{rem}

\section{Simulation Results}
\label{sec:SimResults}
In this section, we validate the proposed controller by simulating IEEE 14-bus system shown in Fig.~\ref{fig:14busTopo}. The circuit specification can be found in~\cite{kodsi2003modeling}. Every generator in IEEE 14-bus example is substituted by VSIs. We assume that the voltage magnitude deviation should be less than 6\% of its nominal value for every bus. In addition, the branch angle difference limit is set to be less than $\gamma = 15\deg$. We find a feasible solution $K$ of Eq.~(\ref{eq: LMIsStability}) by solving Eq.~(\ref{eq:PolynomialBMI}) and maximizing $b$, where $\bar{\mathcal{D}}$ is replaced by $\bar{\mathcal{J}}$ to reduce the complexity as discussed in Remark \ref{rem:ComplexRobust}. The calculated feedback gain $K$ is shown in Table (\ref{table:K14bus}).
\renewcommand{\tabcolsep}{6pt}
\begin{table}[ht]
\caption{Feedback Gains} 
\vspace{2mm}
\centering 
\begin{tabular}{l c c c c c c c} 
\hline\hline 
& &\hspace{-2mm} Unit & \hspace{-3mm} DER 1 & \hspace{-3mm} DER 2 & \hspace{-3mm} DER 3 & \hspace{-3mm} DER 6 & \hspace{-3mm} DER 8 \\ [0.5ex] 
\hline 
\multirow{4}{*}{} & \hspace{-6mm} $K_i(1,1)$ & \hspace{-2mm} mrad/s & \hspace{-3mm} $-5.6$ & \hspace{-3mm} $-1.8$ &\hspace{-3mm} $-10$ &\hspace{-3mm} $-8.8$ &\hspace{-3mm} $-190.7$ \\ 
& \hspace{-6mm} {$K_i(1,2)$} & \hspace{-2mm} mrad/s & \hspace{-3mm} $10.2$ &\hspace{-3mm} $4.9$ & \hspace{-3mm}$8.6$ & \hspace{-3mm}$85$ &\hspace{-3mm} $-24.6$ \\
& \hspace{-6mm} $K_i(2,1)$ & \hspace{-2mm} mV/s &\hspace{-3mm} $-2.7$ &\hspace{-3mm} $-1$ &\hspace{-3mm} $-1.2$ &\hspace{-3mm} $-23.8$ &\hspace{-3mm} $0$ \\ 
& \hspace{-6mm} {$K_i(2,2)$} & \hspace{-2mm} mV/s &\hspace{-3mm} $-0.1$ & \hspace{-3mm}$-1.4$ &\hspace{-3mm} $-10.8$ &\hspace{-3mm} $-16.2$ &\hspace{-3mm} $-40$ \\
\hline 
\end{tabular}
\label{table:K14bus} 
\end{table}
Notice that a linear constraint on $K$ is imposed
\begin{align}
\label{eq:PolyConst}
-b_m \leq K\bar{L}s \leq b_M, \; \forall s \in v(\Xi),
\end{align}
where $\Xi$ represents active and reactive power capacity of inverter buses, $b_M$ and $b_m$ are the upper and lower bounds of ramp of phase angle and voltage. We select $b_M$ and $b_m$ such that $49.7 \leq w \leq 50.3$(Hz) and $|\dot{E}| \leq 0.05$(p.u/s) for all inverters. The constraint on frequency is standard for power networks~\cite{glover2011power}. The upper bound on $|\dot{E}|$ makes the voltage at inverter bus maintain relatively static under any loading condition, while it gradually evolves to the steady state where the power sharing follows. Note that constraint~(\ref{eq:PolyConst}) can be considered as imposing an upper bound of $||K||_2$. Adding constraint~(\ref{eq:PolyConst}) to Eq.~(\ref{eq: LMIsStability}) only changes the optimal value but not the feasibility. As a result, we can always select $b_M$ and $b_m$ only based on physical requirements for the distributed controller design.

\begin{figure}[ht!]
     \begin{center}
        \subfigure[Active power response of the inverters]{%
            \label{fig:P14bus}
            \includegraphics[width=0.5\textwidth]{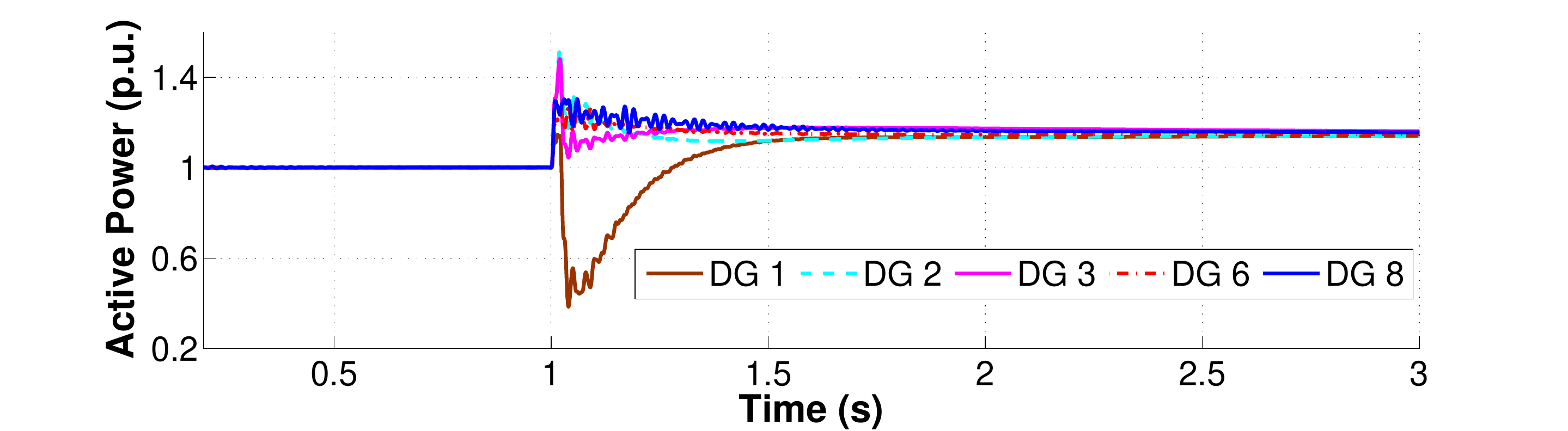}
        }\\%
        \subfigure[Reactive power response of the inverters]{%
           \label{fig:Q14bus}
           \includegraphics[width=0.5\textwidth]{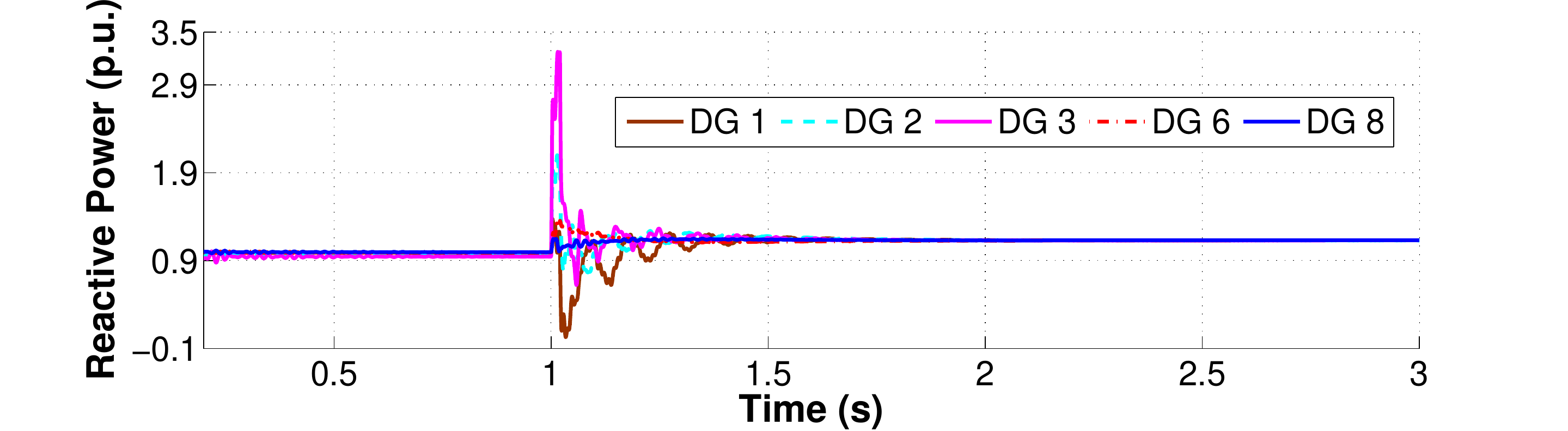}
        }\\ 
        \subfigure[Frequency response of the inverters]{%
            \label{fig:f14bus}
            \includegraphics[width=0.5\textwidth]{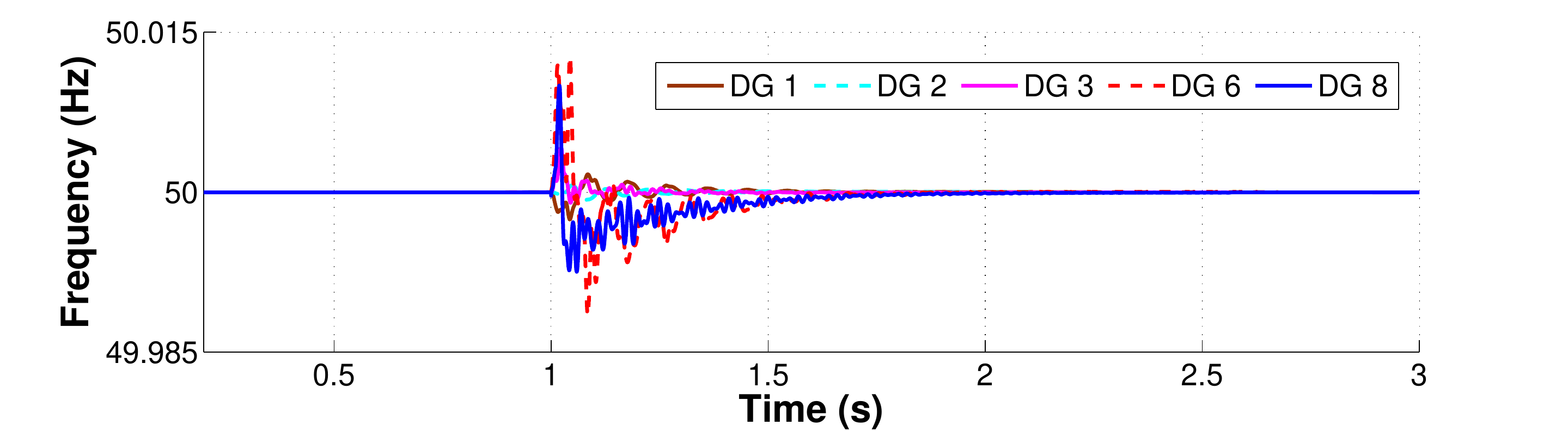}
        }\\%
        \subfigure[Voltage magnitude of buses]{%
            \label{fig:V14bus}
            \includegraphics[width=0.5\textwidth]{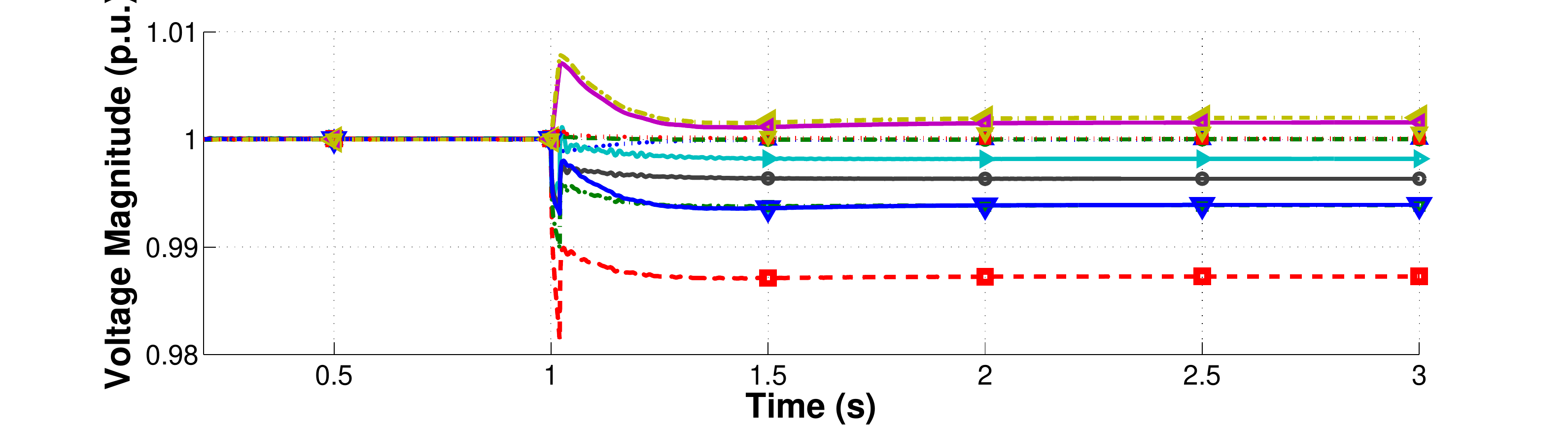}
        }\\%
		\subfigure[Phase angle difference between connected buses]{%
            \label{fig:Ang14bus}
            \includegraphics[width=0.5\textwidth]{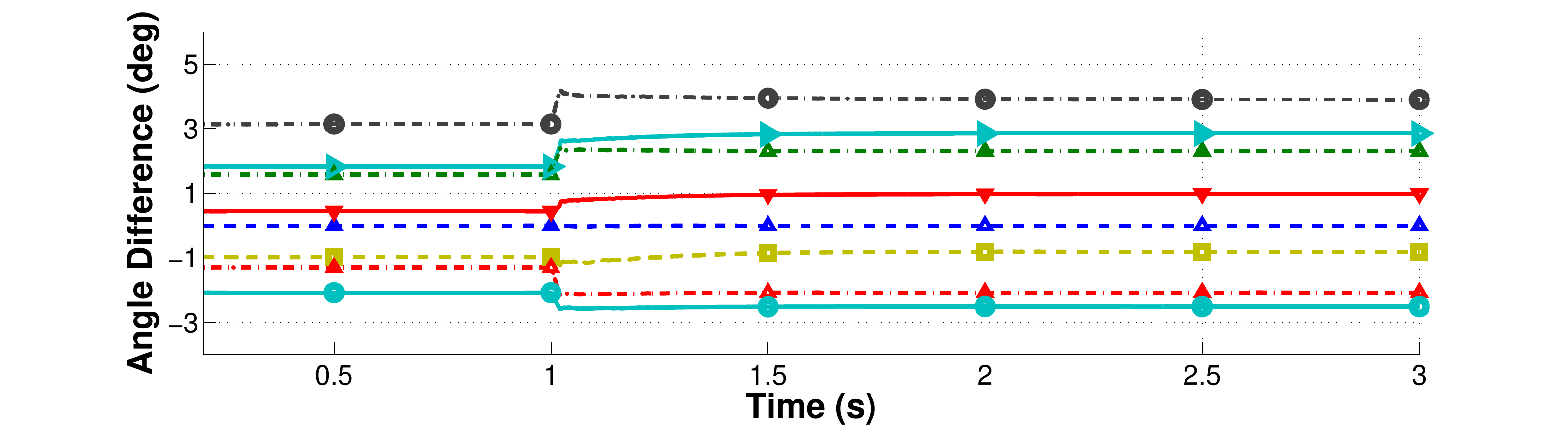}
        }%
    \end{center}
    \caption{%
        Controlled microgrid response when the load at bus 10 changes abruptly at time $t=1$
     }%
   \label{fig:14bus}
\end{figure}

\begin{figure}[ht!]
     \begin{center}
        \subfigure[Active power response of the inverters]{%
            \label{fig:P14bus_LDG1}
            \includegraphics[width=0.5\textwidth]{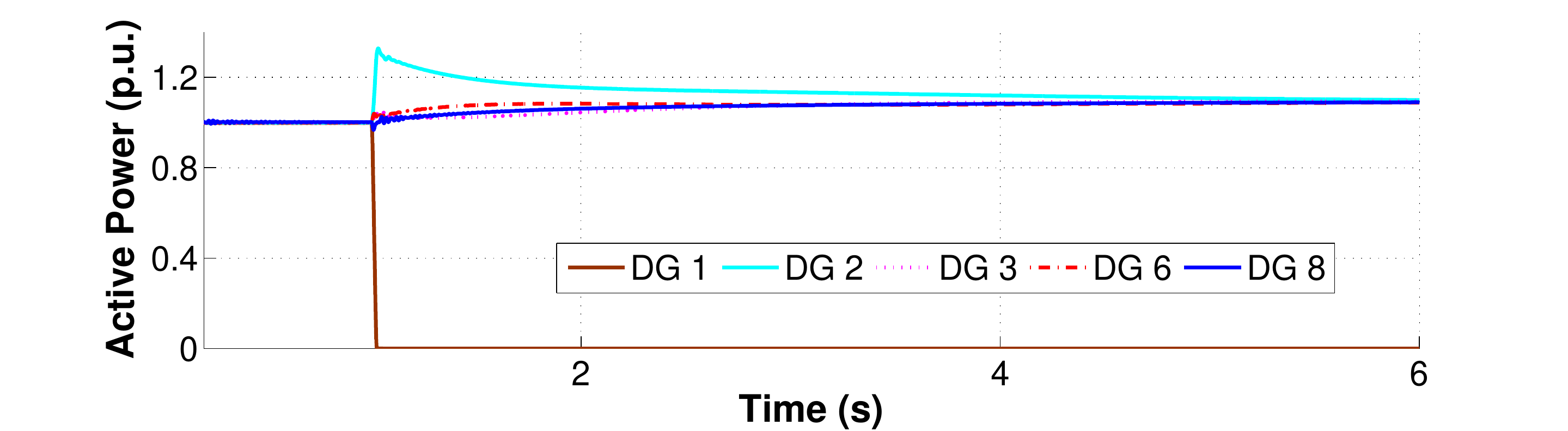}
        }\\%
        \subfigure[Reactive power response of the inverters]{%
           \label{fig:Q14bus_LDG1}
           \includegraphics[width=0.5\textwidth]{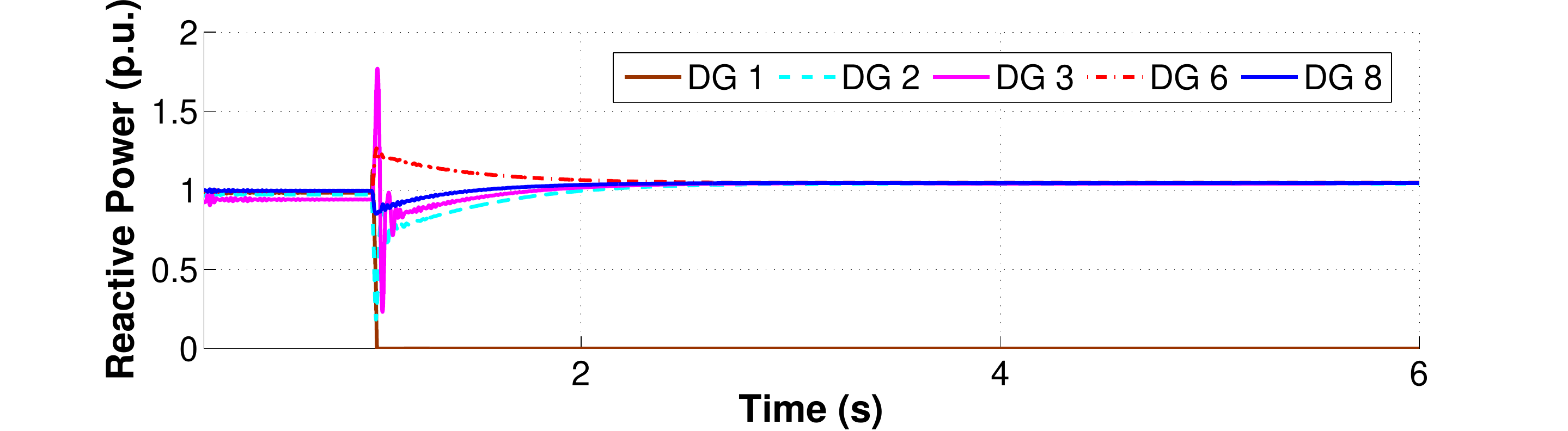}
        }\\ 
        \subfigure[Frequency response of the inverters]{%
            \label{fig:f14bus_LDG1}
            \includegraphics[width=0.5\textwidth]{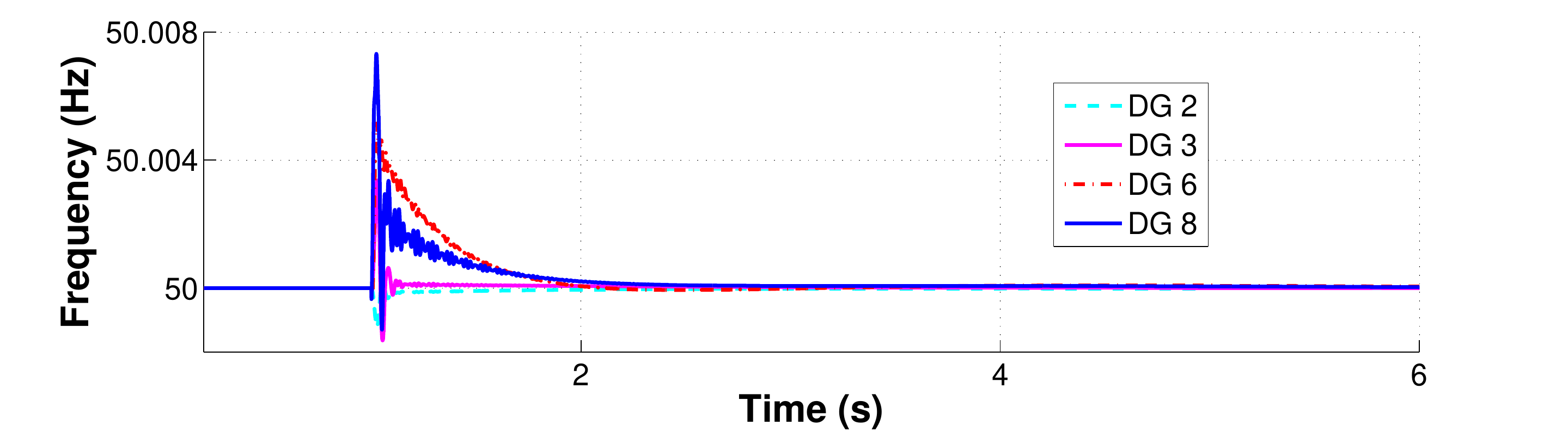}
        }\\%
        \subfigure[Voltage magnitude of buses]{%
            \label{fig:V14bus_LDG1}
            \includegraphics[width=0.5\textwidth]{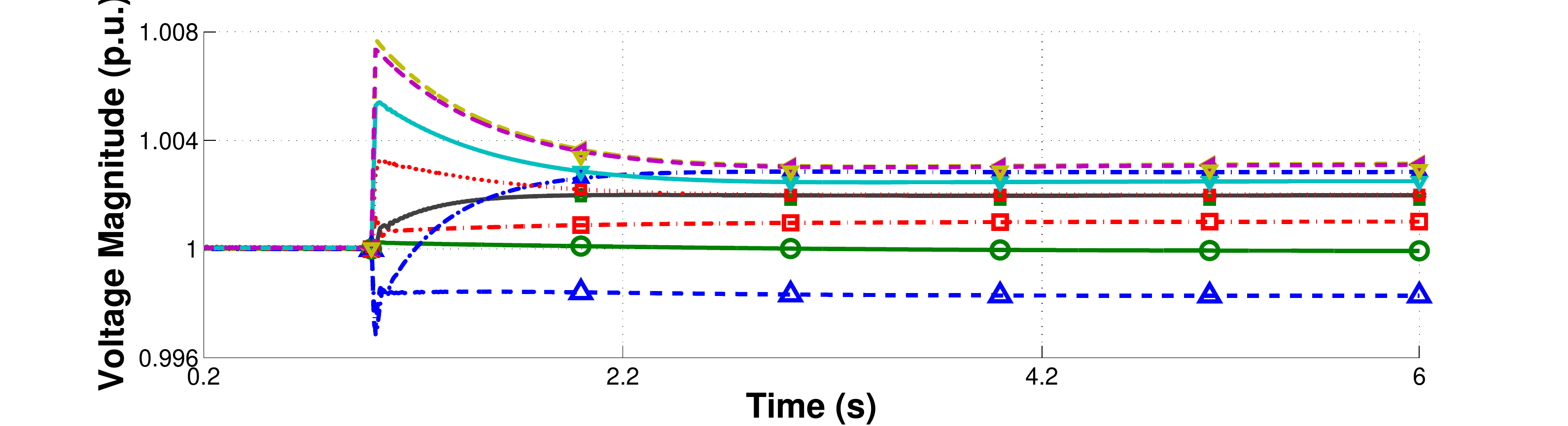}
        }\\%
		\subfigure[Phase angle difference between connected buses]{%
            \label{fig:A14bus_LDG1}
            \includegraphics[width=0.5\textwidth]{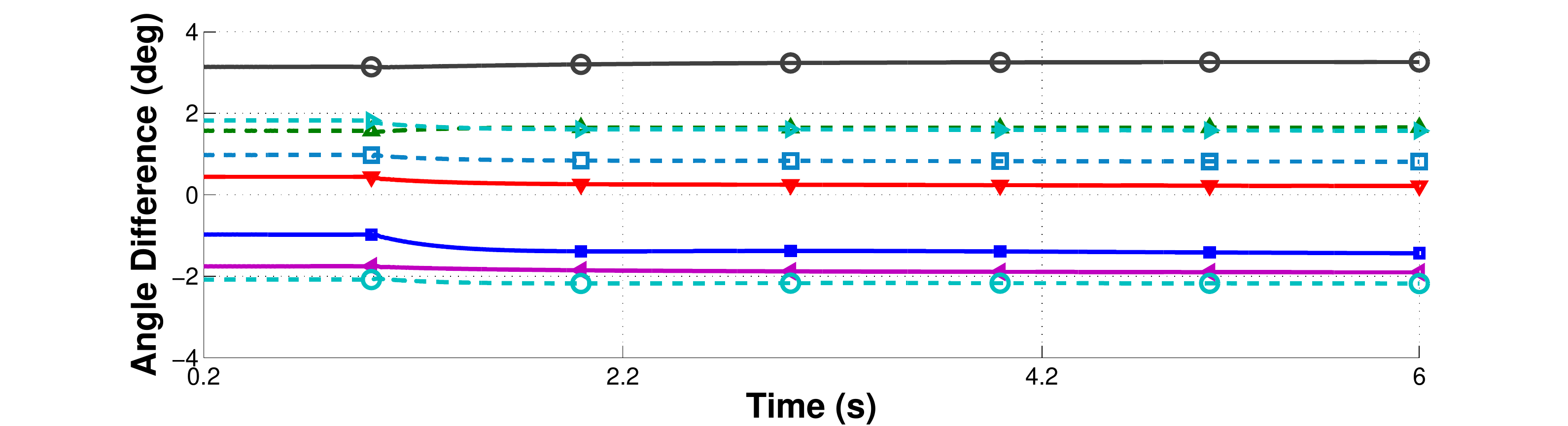}
        }%
    \end{center}
    \caption{%
        Controlled microgrid response when the inverter 1 is failed and disconnected from the microgrid at $1$ seconds.
     }%
   \label{fig:14bus_LDG1}
\end{figure}

We first simulate the case when all the inverters operate normally with an abrupt change of load at bus 10 at time $t=1$. As shown in Fig. \ref{fig:P14bus}- \ref{fig:Ang14bus}, the voltage magnitude always lies inside the desired range. In addition, the mesh microgrid successfully reaches a satisfactory new steady state with proportional power sharing and synchronized frequencies. The second case that we consider is the failure of the inverter 1. From Theorem \ref{thm: Key2}, $K$ solves Eq.~(\ref{eq: LMIsStability}) can maintain the stability and achieve the requirements (\ref{eq:AcPSharing})-(\ref{eq:FReg}) when inverter 1 is disconnected. The simulation results are shown in Fig. \ref{fig:P14bus_LDG1}- \ref{fig:A14bus_LDG1}. It can be seen that the rest of the inverters can carry over the original power injection from the inverter 1 while the desired properties are preserved.

\section{Conclusion}
\label{sec: conclusion}
In this paper, we have presented a distributed control method to coordinate VSIs in an island AC microgrid. Instead of using the conventional droop controller inducing steady state deviation of both frequency and voltage magnitude, our design controls the microgrid to the steady state where frequency synchronization and proportional active and reactive power sharing hold while respecting the voltage regulation constraints. The sufficient conditions of the convergence to the steady state can be approximated into solving optimization problem with BMI constraints. The design is robust with respect to system parameter variations. Unlike most of the existing droop based controller, our design can be applied to the lossy mesh microgrids. Relaxation of the communication requirements and further reduction of the computational complexity are considered as important future works. Extension of the controller for higher order harmonic loads, unbalanced, and other types of nonlinear loads are also among our future works.

\appendix
\section{Proof of Proposition~\ref{assum:Bddxldot}} 
\label{eq:Bddx_L}
We first consider the case of constant impedance loads. Given the impedance at every load bus $i$, the power injection at bus $i$ is determined by local voltage and is written as $P_i^*(x_i)$ and $Q_i^*(x_i)$, for all $i \in \mathcal{V}_L$. Applying KCL at every load bus with power flow equation (\ref{eq:PFEqua}), we have
\begin{align}
\label{eq:KCL_loadbus}
\left\{
\begin{array}{ll}
0 = P_i^*(x_i) + P_i(x) \\ 
0 = Q_i^*(x_i) + Q_i(x)
\end{array}
\right. \text{for all } i \in \mathcal{V}_L.
\end{align}
Taking time derivative of Eq. (\ref{eq:KCL_loadbus})
\begin{align}
\label{eq:JacoAssum}
{\bf{0}}_{2n_L} = f_{I,x}\dot{x}_I + f_{L,x}\dot{x}_L,
\end{align}
where $f_{I,x}$ and $f_{L,x}$ denote the Jacobian matrices follow from Eq. (\ref{eq:KCL_loadbus}). For any $x$, $|| \dot{x}_L ||_2 \leq \kappa_{1,x} || \dot{x}_I ||_2$ always follows from Eq. (\ref{eq:JacoAssum}) if $f_{L,x}$ is invertible. In the case that $f_{L,x}$ is not invertible, a change of coordinate is applied, let $\hat{\dot{x}}_L = T_{L,x}\dot{{x}}_L$, such that $\hat{\dot{x}}_L = [\hat{\dot{x}}_{L_1}^T, \hat{\dot{x}}_{L_2}^T]^T$ and $\hat{\dot{x}}_{L_2}$ span the zero eigenspace of $f_{L,x}$. Let $\kappa_{2,x}$ be the smallest positive eigenvalue of the Hermitian matrix $f_{L,x}^Tf_{L,x}$, then
\begin{align*}
& \kappa_{2,x} \hat{\dot{x}}_{L_1}^T\hat{\dot{x}}_{L_1} \leq \dot{x}_L^Tf_{L,x}^Tf_{L,x}\dot{x}_L  = \dot{x}_I^Tf_{I,x}^Tf_{I,x}\dot{x}_I\\
\Rightarrow &  || \hat{\dot{x}}_{L_1} ||_2 \leq \kappa_{3,x} || \dot{x}_{I} ||_2, \; \kappa_{3,x} \in \mathbb{R}_+.
\end{align*}
Since there exist a simple linear mapping from $\hat{\dot{x}}_{L_1}$ to $\hat{\dot{x}}_{L_2}$, $|| \hat{\dot{x}}_{L_2} ||_2 \leq \kappa_{4,x} || \hat{\dot{x}}_{L_1} ||_2$, we derive
\begin{align*}
||\dot{x}_L ||_2 &= || \hat{\dot{x}}_L ||_2 \leq || \hat{\dot{x}}_{L_1} ||_2 + || \hat{\dot{x}}_{L_2} ||_2 \\
& \leq \kappa_{1,x} || \dot{x}_{I} ||_2,\; \kappa_{1,x} \in \mathbb{R}_+
\end{align*}
for the case that $f_{L,x}$ is not invertible. For $x \in \mathcal{X}_c$, we choose $\kappa = \max_{x\in \mathcal{X}_c}\kappa_{1,x}$ and $|| \dot{x}_L ||_2 \leq \kappa || \dot{x}_I ||_2$ follows. The derivation also applies to constant power loads, only $P_i^*$ and $Q_i^*$ at the load buses become constants instead of the function of local voltage profile, which completes the proof. 

\bibliographystyle{unsrt}      
\bibliography{DroopControl} 
 
\end{document}